\documentclass{article}

  \usepackage[nonatbib,preprint]{neurips_2022}

\usepackage[utf8]{inputenc} 
\usepackage[T1]{fontenc}    
\usepackage{hyperref}       
\usepackage{url}            
\usepackage{array,booktabs}       
\usepackage{amsfonts}       
\usepackage{pifont}
\usepackage{nicefrac}       
\usepackage{microtype}      

\usepackage{natbib}  
\usepackage[dvipsinames]{xcolor,colortbl}
\definecolor{midnightblue}{HTML}{0059b3}
\definecolor{noonblue}{HTML}{e5eef7}
\definecolor{chromered}{HTML}{f14233}
\definecolor{olivedrab}{HTML}{6b8e23}

\usepackage{hyperref} 
            \hypersetup{
                    colorlinks=true,                
                    breaklinks=true,                
                    urlcolor= black,                
                    linkcolor= midnightblue,                
                    bookmarksopen=false,
                    citecolor=black,
                    linkcolor=black,
                    filecolor=black,
                    citecolor=midnightblue,
                    linkbordercolor=white,
}

\usepackage{lipsum}
\usepackage[utf8]{inputenc} 
\usepackage[T1]{fontenc}    
\usepackage{hyperref}       
\usepackage{url}            
\usepackage{booktabs}       
\usepackage{amsfonts}       
\usepackage{nicefrac}       
\usepackage{microtype}      
\usepackage{xcolor}
\usepackage{graphicx}
\usepackage{subfigure}
\usepackage{booktabs} 
\usepackage{hyperref}
\usepackage{times}    
\usepackage[T1]{fontenc}
\usepackage[utf8]{inputenc} 
\usepackage{amssymb,amsmath,amscd,amsfonts,amsthm,bbm,mathrsfs,yhmath}

\usepackage{caption}
\usepackage{graphics}
\usepackage{mathtools}
\usepackage{cleveref}
\usepackage{comment}

\usepackage{paralist,enumitem}                                       
\usepackage{pdfpages}

\usepackage[titletoc,toc,page]{appendix}
\usepackage{minitoc}

\newtheorem{theorem}{Theorem}
\newtheorem{lemma}{Lemma}
\newtheorem{corollary}{Corollary}
\newtheorem{proposition}{Proposition}
\newtheorem{definition}{Definition}
\newtheorem{remark}{Remark}

\newtheorem{assumption}{Assumption}

\newtheorem{assumptionalt}{Assumption}[assumption]
\newenvironment{assumptionp}[1]{
  \renewcommand\theassumptionalt{#1}
  \assumptionalt
}{\endassumptionalt}

\usepackage[colorinlistoftodos,bordercolor=orange,backgroundcolor=orange!20,linecolor=orange,textsize=scriptsize]{todonotes}

\newcommand{\rmd}{{\mathrm d}}
\newcommand{\rmP}{{\mathrm P}}
\newcommand{\rmJ}{{\mathrm J}}
\newcommand{\rmT}{{\mathrm T}}
\newcommand{\rmQ}{{\mathrm Q}}
\newcommand{\cF}{{\mathcal F}}

\newcommand{\cO}{{\mathcal O}}

\newcommand{\cN}{{\mathcal N}}

\newcommand{\cP}{{\mathcal P}}

\newcommand{\cH}{{\mathcal H}}

\newcommand{\R}{{\mathbb R}}

\newcommand{\KL}{\mathop{\mathrm{KL}}\nolimits}

\newcommand{\Stein}{\rm Stein}

\newcommand\NN{\mathbb N}
\newcommand\RR{\mathbb R}
\newcommand\EE{\mathbb E}

\newcommand{\norm}[1]{\left\| #1 \right\|}

\newcommand{\inner}[2]{\left< #1 , #2 \right>}

\title{Convergence of Stein Variational Gradient Descent under a Weaker Smoothness Condition}

\author{%
  Lukang Sun \\
  CEMSE, KAUST \\
  \texttt{lukang.sun@kaust.edu.sa}\\
  \And 
  Avetik Karagulyan \\
  CEMSE, KAUST \\
  \texttt{avetik.karagulyan@kaust.edu.sa}\\
  \And
   Peter Richtarik \\
  CEMSE, KAUST \\
  \texttt{peter.richtarik@kaust.edu.sa}\\
}

\begin{document}

\maketitle

\begin{abstract}
  Stein Variational Gradient Descent (SVGD) is an important alternative to the Langevin-type algorithms for sampling from probability distributions of the form 
  $\pi(x) \propto \exp(-V(x))$. 
  In the existing theory of Langevin-type algorithms and SVGD, the potential function $V$ is often assumed to be $L$-smooth.
  However, this restrictive condition excludes a large class of potential functions such as polynomials of degree greater than $2$. 
  Our paper studies the convergence of the SVGD algorithm for distributions with $(L_0,L_1)$-smooth potentials. 
  This relaxed smoothness assumption was introduced by \cite{zhang2019gradient} for the analysis of gradient clipping algorithms. With the help of trajectory-independent auxiliary conditions, we provide a descent lemma establishing that the algorithm decreases the $\KL$
  divergence at each iteration and prove a complexity bound for SVGD in the population limit in terms of the Stein Fisher information.

\end{abstract}
\addtocontents{toc}{\protect\setcounter{tocdepth}{0}}
\section{Introduction}\label{sec:intro}
 Bayesian  methods are widely implemented in various inference tasks that emerge in computational statistics 
  \citep{neal1992bayesian,roberts1996exponential}, 
  machine learning \citep{grenander1994representations,csimcsekli2017fractional}, 
  inverse problems \citep{zhou2020bayesian,durmus2018efficient} and model selection \citep{feroz2013exploring,cai2021high}.
  Often Bayesian methods boil down to approximate integration problems which are solved using Markov Chain Monte-Carlo algorithms \citep{robert1999monte}. 
  In practice, the target distribution  $\pi$ defined on $\RR^d$
  is often absolutely continuous w.r.t. the Lebesgue measure and we have access to its 
  density up to a normalization constant.  
  Throughout the paper we will use the same notation for probability distributions and their corresponding densities. In Bayesian statistics the target distribution is generally given by 
  \begin{equation*}
    \label{eq:target}
    \pi(x)=\frac{1}{Z} e^{-V(x)},
  \end{equation*}
  where $Z$ is the normalization constant and $V:\R^d\to \RR$ is called potential function. Here, we usually do not know the value of $Z$. 
  The approximate sampling methods are proceed as follows. 
  For every $\varepsilon > 0$, the goal is to construct a distribution $\mu$ that approximates the target $\pi$ in some probability distance {\sf dist}: 
  \begin{equation*}
    {\sf dist}(\mu,\pi) < \varepsilon.
  \end{equation*}
  In this paper, we will study the case when ${\sf dist}$ is the Kullback-Leibler divergence:
  \begin{equation*}
    \label{eq:KL}
    \KL(\mu\mid\pi):=\int \log\left(\frac{\mu(x)}{\pi(x)}\right)\mu(\rmd x) 
    = \EE_{X\sim\mu}\big[V(X)\big] - H(\mu).
  \end{equation*}
  Here, $H(\mu) := - \int \log(\mu(x)) \mu(\rmd x)$ is the negative entropy. 
  The sampling problem can also be seen as a minimization problem in the space of probability measures (see e.g. \citep{liu2016stein, wibisono2018sampling,durmus2018analysis}). 
  Indeed, consider the functional $\cF(\cdot) := \KL (\cdot \mid \pi)$. 
  This functional is non-negative and it admits its minimum value $0$ only when for the target distribution:
  $\pi = \arg\min_{\mu} \cF(\mu)$. 
  The connection of sampling and optimization has been repeatedly leveraged in the previous literature of Langevin sampling.
  Various algorithms such as Langevin Monte-Carlo (e.g. \cite{dalalyan2017further,wibisono2019proximal,durmus2018analysis}), Underdamped Langevin Algorithm (e.g. \cite{ma2019there,chatterji2018theory})  have often  been influenced by the known optimization methods. 
  Conversely, another line of research has studied the application of sampling algorithms to solve optimization problems (see e.g. \citep{raginsky2017non}).
  In this paper, we will study another algorithm called Stein Variational Gradient Descent (SVGD), which is designed as a gradient descent algorithm for the $\KL$ divergence in the space of probability measures.
  {}

\subsection{SVGD}
  
  The LMC algorithm, treats the functional $\KL(\cdot|\pi)$ as a composite functional described in \eqref{eq:KL}. 
  For details, we refer the reader to \Cref{sec:app-langevin} (see also \citep{wibisono2018sampling}).
  Unlike the LMC, the Stein Variational Gradient Descent (SVGD)
  algorithm applies gradient descent directly to $\KL(\cdot|\pi)$ (see \Cref{sec:SVGD} for the complete definition). 
  SVGD is an important alternative to the Langevin algorithm and already has been used extensively  in different settings of  
  machine learning, such as variational auto-encoders \citep{pu2017vae}, reinforcement learning~\citep{liu2017policy}, 
  sequential decision making~\citep{zhang2018learning,zhang2019scalable}, generative adversarial networks~\citep{tao2019variational} and federated learning~\citep{kassab2020federated}. 

  The seminal work of \cite{liu2016stein} introduced SVGD as a sampling method. 
  Since then several variants of SVGD have been proposed. 
  Here is a non-exclusive list: random batch method SVGD \citep{li2020stochastic}, matrix kernel SVGD~\citep{wang2019stein}, Newton version SVGD \citep{detommaso2018stein}.
  However, the theoretical understanding of SVGD is still limited to the continuous time approximation of SVGD or the so called infinite particle regime.
  The work by \cite{liu2017stein} was the first that proved a convergence result of the SVGD in the population limit. 
  Later, \cite{duncan2019geometry} studied the geometry related to the SVGD and proposed a scheme to choose the kernel.  
  The clean analysis of the finite particle regime remains an open question (see also \citep{lu2019scaling,nusken2021stein}). 
  In \citep{korba2020non}, a descent lemma was established  for the SVGD in population limit in terms of Kullback-Leibler 
  divergence. 
  The drawback of this result was that the analysis relied on the path information of the SVGD which is unknown 
  beforehand. 
  \cite{salim2021complexity} improved the work of ~\cite{korba2020non} and provided a clean analysis for the convergence. 
  They assumed $\pi$ satisfies $\text{T}_1$ (see \Cref{sec:setting}) inequality which essentially replaced the initial trajectory condition. 
  This implied a complexity bound for the SVGD in terms of the desired 
  accuracy $\epsilon$ and dimension $d$.  To the best of our knowledge, this is the state 
  of the art result on SVGD.


\subsection{Contributions}
  The main contribution of the paper relies on its weaker set of assumptions, that allow to treat a larger class of probability distributions which includes densities with polynomials. 
  We enlarge the class of probability distributions two-fold. 
  \paragraph{Smoothness:} 
  The gradient smoothness assumption is very common in the sampling literature (see e.g., 
  \cite{durmus2017nonasymptotic, dalalyan2017theoretical,
  dalalyan2019user, durmus2018analysis, vempala2019rapid, shen2019randomized,korba2020non,salim2021complexity}).
  It is formulated as the Lipschitz continuity of the gradient. That is
  the Hessian $\nabla^2 V(\cdot)$ is well-defined on $\RR^d$ and 
  \begin{equation}\label{eq:L-smooth}
    \left\| {\nabla^2 V}(x) \right\|_{o p} \leq L, \quad \forall x\in\mathbb{R}^d.
  \end{equation}     

  Though people have made great progress towards the understanding of these algorithms, the $L$-smoothness condition is quite restrictive. 
  In fact from $L$-smoothness condition we can easily deduce that $V$ has at most quadratic growth rate. 
  In particular, the large class of polynomial potentials whose order is higher than 2 does not satisfy this condition. 
  For the LMC algorithm, several papers have proposed different methods to relax or remove this assumption \citep{lehec2021langevin,brosse2019tamed,hutzenthaler2012strong,sabanis2013note,erdogdu2021convergence,chewi2021analysis}. 


  However, to the best of our knowledge there is no result for the SVGD under a relaxed smoothness assumption. 
  In this paper we introduce the concept of $(L_0,L_1)$ smoothness (see \Cref{A-l0l1} for the definition). The latter was initially proposed by  \cite{zhang2019gradient}, for the gradient clipping algorithm. 
  The new condition with parameters $(L,0)$ is equivalent to the $L$-smoothness and therefore, it is indeed a weaker assumption. 
  Important example of such functions are the higher order polynomials (see \Cref{sec:app-applications} for details).

  \paragraph{Functional inequalities:}

  The analysis of Langevin algorithms often relies on the strong convexity of the
  potential function $V$. 
  A line of work has proposed different methods to relax or bypass this assumption. 
  One possible solution is to slightly modify the original algorithm (e.g. \citep{dalalyan2019bounding,karagulyan2020penalized}). 
  Another approach relies on functional inequalities such as Poincaré or logarithmic Sobolev inequalities (e.g. \citep{vempala2019rapid,chewi2020exponential,ahn2021efficient}).
  However,  the verification of these inequalities is not straightforward. 

  Unlike the Langevin algorithm, the SVGD algorithm does not require convexity in any form. 
  However, functional inequalities also serve for the analysis of the SVGD but with a different purpose. 
  \cite{salim2021complexity} have proved a complexity result in dimension and precision error for the SVGD algorithm
  performed on the targets that satisfy Talagrand's $\rmT_1$ inequality (see \eqref{eq:Tp} for details). 
  This assumption  replaces the bound on the trajectory that was initially proposed by \cite{korba2020non}. 
  Despite the major improvement, they still assumed the $L$-smoothness of the potential, which as mentioned previously does not cover the polynomials.

  In this regard, we propose to replace the Talagrand's inequality with the generalized $\rmT_p$ (see \Cref{Tp-assumption} for the definition). Due its general form, this new condition includes inequalities that are easy to verify (see \Cref{sec:app-assumption}). 
  See also the \Cref{Table:contribution} for a visual representation of the literature review.

  \begin{table}[ht]\scriptsize 
  \centering
  \renewcommand{\arraystretch}{1.5}
    \begin{tabular}{>{\centering}m{.5in} >{\centering}m{.25in} >{\centering}m{.6in}  > {\centering}m{.5in} > {\centering}m{.55in} > {\centering}m{1.25in} > {\centering\arraybackslash}m{.3in}}
      \toprule
      Paper & Method & Smoothness  & Other  conditions & Can deal with higher-order polynomials?& Complexity to get $\varepsilon$ error & Criterion \\
      \bottomrule
      \cite{vempala2019rapid} & LMC & $L$-smoothness  & $\lambda_{LS}$-Log-Sobolev & {\color{chromered}\ding{55}} &
             $\widetilde{\cO}\left[\frac{d}{\lambda_{LS}^2 \varepsilon}\right]$ & $\KL$\\
      \hline
      \cite{chatterji-lmc-smoothness} & LMC & $(L,\alpha)$-Hölder continuity + Gaussian smoothing & Log-concavity & {\color{chromered}\ding{55}} &  $\widetilde{\mathcal{O}}\left[\frac{d^{\nicefrac{(5-3\alpha)}{2}}}{\varepsilon^{\nicefrac{4}{(1+\alpha)}}}\right]$ & $W_2$ \\
       \hline
       \cite{chewi2021analysis} & LMC & Hölder  continuity &  $\alpha$-tails + Modified  Log-Sobolev & {\color{chromered}\ding{55}} 
       & $\widetilde{O}\left[\frac{d^{(2 / \alpha)(1+1 / s)-1 / s }}{\varepsilon^{1 / s}}\right] $ & $\KL$ \\
      \bottomrule
       \cite{korba2020non} & SVGD & $L$-smoothness & Trajectory bounds &  {\color{chromered}\ding{55}} & 
       $\cO\left[ \frac{d}{\varepsilon} \sqrt{C} \right]$ & $I_{\Stein}$\\
       \hline
      \cite{salim2021complexity} & SVGD & $L$-smoothness & $\rmT_1$ inequality with constant $\lambda_{\rmT}$ & {\color{chromered}\ding{55}}
       & $\widetilde{\mathcal{O}}\left[\frac{d^{3 / 2}}{\lambda_{\rmT}^{1 / 2} \varepsilon}\right]$ &$I_{\Stein}$\\
      \hline
      \rowcolor{noonblue}  This paper & SVGD & $(L_0,L_1)$ smoothness &   $\rmT_p$ inequality with constant $\lambda_{\rmT}$
      &  {\color{olivedrab}\ding{51}} &$\cO\Big[ {(\lambda_{\rmT}^{\frac{p}{2}} \varepsilon)}^{-1}{ (p d)^{\frac{(p+1)(p+2)}{4}}}\Big] $ & $I_{\Stein}$\\
      \hline
      \rowcolor{noonblue}  This paper & SVGD & $(L_0,L_1)$ smoothness & Generalized \ref{Tp-assumption} 
      &  {\color{olivedrab}\ding{51}} & $\cO \left[\frac{\lambda_{BV}^p (pd)^{p+1}}{\varepsilon}\right]$  & $I_{\Stein}$\\
      \hline
    \end{tabular}
    \vspace{.3cm}
    \caption{This summarizes several previous results on the LMC and SVGD algorithms.
    Hear $\tilde{\cO}$ corresponds to the order without the log-polynomial factors.  
    One has to bear in mind, that the complexities of the SVGD and LMC
    cannot be compared to each other as they use different error metrics. 
    The complexity of \citep{korba2020non}[Corollary 6] contains the trajectory bound constant $C$, which may depend on the dimension $d$.
    In the last row of the table, the \ref{Tp-assumption} is satisfied for the function $\rmJ(r) = \lambda_{BV}(r^{1/p} + (r/2)^{1/2p})$, where $\lambda_{BV}$ is a constant that may depend on the dimension (see \Cref{sec:assumptions}).}
    \label{Table:contribution}
  \end{table}

\subsection{Paper structure}
  
  The paper is organized as follows. 
  In \Cref{sec:setting}, we present the  mathematical setting of our problem. 
  We  introduce the basic notions, describe the main algorithm and list the necessary assumptions. 
  The first part of \Cref{sec:main} is devoted to the main result, which is a descent lemma for $\KL(\cdot,\pi)$ functional on the Wasserstein space. 
  The second part contains the complexity results and discussions.
  \Cref{sec:app-applications} presents two examples.
  \Cref{sec:conclusion} concludes the results of the paper and discusses possible future work.
  Finally in the Appendix, the reader may find the postponed proofs, intermediary technical lemmas and supplementary discussions.

\section{Mathematical setting of the problem}\label{sec:setting}

\subsection{Notations}

  We will denote the random variables with uppercase (resp. lowercase) Latin letters the random (resp. deterministic) 
  vectors, unless specified. The space of $d$-dimensional real vectors is denoted by $\RR^d$, while the set of non-negative real numbers is denoted
  by $\RR^+$.
  Our target distribution $\pi$ is defined on $\RR^d$, which is enhanced with the Euclidean norm. 
  The notation $\|\cdot\|$ will correspond to the ${\ell}_2$ norm on $\RR^d$ unless specified. 
  Also, we will assume that $\pi$ has the $p$-th moment, that is  $\pi \in \mathcal{P}_{p} (\RR^d)$.   
  The Jacobian of a vector valued function 
  $h(\cdot) = (h_1(\cdot),\ldots,h_d(\cdot))^\top$ is a $d\times d$ matrix defined as
  \begin{equation*}
    J h(x) := \big({ \partial_{x_i} h_j}\big)_{i=1,j=1}^{d,d}.
  \end{equation*}
  The divergence of the vector valued function $h$ is the trace of its Jacobian:
  \begin{equation*}
    \operatorname{div} h(x) := \sum_{i=1}^{d}  \partial_{x_i} h_i (x).
  \end{equation*}
  The Hessian of a real valued function $U : \RR^d \rightarrow \RR$ is defined as the following square matrix:
  \begin{equation*}
    \nabla^2 U (x) := \left(\frac{\partial^2 h}{\partial_{x_i} \partial_{x_j}} (x) \right)_{i=1,j=1}^{d,d}.
  \end{equation*}

  For any Hilbert space $\mathcal{H}$, we denote by $\langle\cdot, 
  \cdot\rangle_{\mathcal{H}}$ the inner product of $\mathcal{H}$ and by $\|\cdot\|_{\mathcal{H}}$ its norm. Moreover, $\|\cdot\|_{\mathrm{op}}$ denotes the operator norm on the set of matrices. 
  Let  $\mu,\nu \in \mathcal{P}_{p} (\RR^d)$. The set of couplings $\Gamma(\mu,\nu)$ is the set of all 
  joint distributions defined on $\RR^d \times \RR^d$ having $\mu$ and $\nu$ as its marginals.
  The Wasserstein-$p$ distance between two probability measures is defined as 
  \begin{equation*}
    W_p(\mu,\nu) := \inf _{\eta \in {\Gamma}(\mu, \nu)} \left[\int\|x-y\|^{p} \eta(\rmd x, \rmd y)\right]^{1/p}.
  \end{equation*}
  The Kullback-Leibler divergence is defined as
  \begin{equation*}
    \KL (\mu\mid\nu) =
    \begin{cases}
       \int_{\RR^d}  \log\left(\frac{\mu(x)}{\nu(x)}\right)\mu(\rmd x), \text{ if } \mu \ll \nu;\\
       +\infty, \text{ otherwise}.
    \end{cases}
  \end{equation*}
  We will use  the spectral and the Hilbert-Schmidt norms for matrices. For $M \in \RR^{d\times d}$ they are respectively defined as
  \begin{equation*}
    \begin{aligned}
      \|M\|_{op} &:= \sqrt{\lambda_{\max} (M^{\top}M)}, \qquad \text{and} \qquad
      \|M\|_{HS} &:= \sqrt{\sum_{i=1}^{d} \sum_{j=1}^{d} M_{i,j}^2}.
    \end{aligned}
  \end{equation*}
  Here $\lambda_{\max}$ corresponds to the largest eigenvalue.


\subsection{The definition of the SVGD}\label{sec:SVGD}

  Let us present briefly  Reproducing Kernel Hilbert Spaces (RKHS) and some of its essential properties. 
  We refer the reader to \citep{steinwart2008support}[Chapter 4] for a detailed introduction. 
  Let the map $k : \RR^d \times \RR^d \rightarrow \RR$ be a reproducing kernel and let $\cH_0$ be its corresponding RKHS. 
  This means that $\cH_0$ consists of real-valued maps from $\RR^d$ to $\RR$, including the feature maps $\Phi(x) := k(x,\cdot) \in \cH_0$, and the reproducing property is satisfied:
      \begin{equation*}
        f(x) = \langle f, k(\cdot,x) \rangle_{\cH_0}.
      \end{equation*}
  See the discussion on the reproducing property in \Cref{append:repro}.
  Let $\cH$ be the space of the $d$-dimensional maps $\{(f_1,\ldots,f_d)^{\top} \mid f_i \in \cH_0, i = 1,\ldots,d\}$.
  For two vector functions $f = (f_1,\ldots,f_d)^{\top}$ and $g = (g_1,\ldots,g_d)^{\top}$ from $\cH$, we define the scalar product as
  \begin{equation*}
    \langle f,g \rangle_{\cH} := \sum_{i=1}^{d} \langle f_i,g_i \rangle_{\cH_0}.
  \end{equation*}

  Suppose that we have a kernel $k : \RR^d\times\RR^d \rightarrow \RR^+$ and $\cH_0$ 
  is its corresponding RKHS. As described above, we construct the Hilbert space $\cH = \cH_0^d$. 
  Our goal is to construct an iterative algorithm, that has its iterates in the space of probability measures 
  $\cP_p(\RR^d)$. 
  Each iterate is  defined as a pushforward measure from the previous one in a way that it minimizes 
  the $\KL$ distance the most.
  To do so, for every $\psi \in \cH$ and $\gamma > 0$ let us define the operator  
  \begin{equation*}
    T_{\gamma}(x) := x - \gamma \psi(x).
  \end{equation*}
  The operator $\psi$ will serve us as the direction or the perturbation, while as $\gamma$ is the step-size. 
  The goal is to choose the direction in which the $\KL$ error descends the most. 
  Thus, for every $\mu$, the optimal choice of $\psi$ is the solution of the following problem:
  \begin{equation}\label{eq:opt-dir-problem}
    g_\mu := \arg\max_{\psi \in \mathcal{H}}\left\{- \frac{\mathrm{d}}{\mathrm{d} \gamma} \mathrm{KL}\left(T_{\gamma}{\#} \mu \| \pi \right)\Big|_{\gamma=0} \quad \text { s.t. } \quad\|\psi\|_{\mathcal{H}} \leq 1\right\} .
  \end{equation}
  The functional objective depends linearly on the perturbation function $\psi$. Indeed,
  \cite{liu2016stein} have proved that 
  \begin{equation}\label{eq:linear-obj}
    - \frac{\mathrm{d}}{\mathrm{d} \gamma} \mathrm{KL}\left(T_{\gamma}{\#} \mu \| \pi \right)\Big|_{\gamma=0}
    = \int_{\RR^d}  \big(- V(x)\psi(x) + \operatorname{div} \psi(x)\big) \mu(\rmd x).
  \end{equation}
  If $\psi$ satisfies mild conditions, then the right-hand side of \eqref{eq:linear-obj}  is equal to zero if and only if the measures $\mu$ and $\pi$ coincide 
  \citep{stein1972bound}. 
  This property motivates to define a discrepancy measure called Stein discrepancy as follows:
  \begin{equation*}
    I_{\Stein}(\mu \mid \pi) := \max_{\psi \in \mathcal{H}} 
      \left\{ 
      \int_{\RR^d}  \big(- V(x)\psi(x) + \operatorname{div} \psi(x)\big) \mu(\rmd x) \quad \text { s.t. } \quad\|\psi\|_{\mathcal{H}} \leq 1
      \right\}^2.
  \end{equation*}
  \cite{liu2016kernelized} have proved that the solution of the optimization problem \eqref{eq:opt-dir-problem}  is  given explicitly by
  \begin{equation*}
    g_\mu (\cdot) = - \int_{\RR^d} \big[\nabla\log \pi(x) k(x,\cdot)  + \nabla_x k(x,\cdot)\big]\mu (\rmd x)
  \end{equation*}
  and  $\|g_{\mu}\|_{\cH} = \sqrt{I_{\Stein} (\mu \mid \pi)}$. 
  For the proof of this inequality, we refer the reader to the remark in \Cref{sec:propproof}. In our case,  $I_{\Stein}$ is often referred to as the squared Kernelized Stein Discrepancy (KSD) or the Stein-Fisher information. 
  Under certain conditions (see \Cref{sec:app-weak-conv}) on the target distribution $\pi$ and the kernel $k$, the convergence in KSD implies weak convergence:
  \begin{equation*}
      I_{\Stein} (\mu_n \mid \pi) \rightarrow 0  \implies \mu_n \Rightarrow \pi.
  \end{equation*}

  \begin{remark}\label{rem:alt-def-gmu}
    Integration by parts yields the following formula for $g_{\mu}$:
    \begin{equation*}
      g_{\mu} (\cdot) = \int_{\RR^d} \nabla \log\left(\frac{\mu(x)}{\pi(x)}\right) k(x,\cdot) \mu(\rmd x).
    \end{equation*}
  \end{remark}
  This latter equality is the alternative definition of the optimal direction $g_\mu$ given by \cite{korba2020non}. 
  The proof of the remark can be found in \Cref{sec:app-proof-alt-def-gmu}.
  Thus, we have determined the descent direction. In order to formulate the algorithm let us initialize our algorithm
  at $\mu_0 \in \cP_{p}(\RR^d)$. The $(i+1)$-th iterate $\mu_{i+1}$ is obtained by transferring $\mu_i$ with the map
  $I - \gamma g_{\mu_i}$:
  \begin{equation}\label{eq:SVGD}\tag{SVGD}
    \mu_{i+1} = (I - \gamma g_{\mu_i})\# \mu_i.
  \end{equation}
  In case when $k \in L_p(\mu)$, then $\cH_0 \subset L_p(\mu)$ (see \citep{steinwart2008support}[Theorem 4.26]). 
  This guarantees that every iteration of the SVGD has a finite $p$-th moment, that is for every $i\in \NN$,  $\mu_i \in \cP_p(\RR^d)$.
  The rest of the section is devoted to the assumptions that we will use later in the analysis. 

%
  
\subsection{Assumptions}\label{sec:assumptions}
  In order to perform SVGD  we need the target $\pi$, in particular the gradient of its potential potential $V$ and the kernel 
  $k(\cdot,\cdot)$.   
  Below we present the four main assumptions on $\pi$,$V$ and $k$, that we use in the analysis. 
  The first assumptions is a relaxed form of smoothness condition for $\nabla V$. 
  
  We study the SVGD algorithm under the so called $(L_0,L_1)$ smoothness, which has been borrowed from optimization literature.
  In \citep{zhang2019gradient}, the authors studied the convergence of the clipped gradient descent under this condition. 
  The assumption goes as follows.
  \begin{assumptionp}{$(L_0,L_1)$}\label{A-l0l1}
    The Hessian ${\nabla^2 V}$ of $V=-\log \pi$ 
    is well-defined and $\exists L_0, L_1\geq 0$ s.t. 
      \begin{equation}\label{eq:l0l1}
        \left\| {\nabla^2 V}(x) \right\|_{op} \leq L_0+L_1 \| \nabla V(x) \|,
      \end{equation}
     for any $x\in\mathbb{R}^n$.
  \end{assumptionp}
  First let us notice that when $L_1 = 0$ the inequality \eqref{eq:l0l1} becomes the usual smoothness condition \eqref{eq:L-smooth}. 
  On the other hand, let us consider the generalized normal distribution (also known as exponential power distribution) which is given by its density 
  $\pi(x) \propto \exp(-\|x-a\|^{\beta})$, where  $\beta \geq 2$. 
  In this case, Hessian $\nabla^2 V$ is not bounded by a constant in terms of the operator norm. 
  However, it is easy to verify that $V$ satisfies \Cref{A-l0l1}, which yields that this assumption is indeed a relaxation of the $L$-smoothness.

  The next two assumptions are intended to  replace the trajectory condition \citep{korba2020non,duncan2019geometry} for  polynomial potentials.
  However, this condition limits the applicability of the algorithm as these results cannot guarantee the convergence 
  of the algorithm before actually implementing it. 
  This condition was replaced by Talagrand's $\rmT_1$ inequality in the recent work by \cite{salim2021complexity}, 
  where complexity bounds are proved for $L$-smooth potentials.  
  Talagrand's $\rmT_p$ inequality goes as follows: 
  \begin{equation}\label{eq:Tp}
    W_p(\mu,\pi)\leq \sqrt{\frac{2\KL(\mu\mid\pi)}{\lambda_{\rmT}}}, \quad \forall \mu\in\cP_p(\R^d).
  \end{equation}
  Since the Wasserstein-$p$ distance is increasing with respect to the order $p$, then $\rmT_p$ implies $\rmT_1$. 
  In general, it is hard to verify if the distribution satisfies the $\rmT_p$ inequality. 
  However, one may check that $\rmT_2$ is true for $m$-strongly log-concave distributions with $\lambda_{\rmT} = m$.
  According to Bakry-Émery theory,  sharp estimates on this constant are available due to perturbation arguments, such as the well-known Holley-Stroock method 
  (see \cite{steiner2021feynman}).
  We propose a more general form of the $\rmT_p$ inequalities which will help us to analyze the convergence of the SVGD under the $(L_0,L_1)$ smoothness condition.
  \begin{definition}[Generalized $\rmT_p$ inequality]
    Let $p \geq 1$. The distribution $\pi$ satisfies the  generalized $\rmT_p$ inequality if there exists  
    an increasing function $\rmJ : \RR^+ \rightarrow \RR^+$
    such that for all $\mu \in \mathcal{P}_{p}(\RR^d)$, we have $W_{p}(\mu, \pi) \leq \rmJ(\rm{KL}(\mu|\pi))$, where $W_p$ is the Wasserstein-$p$ distance.  
  \end{definition}

  \begin{assumptionp}{$(\rmT_p,{\rm J})$}\label{Tp-assumption}
    The target distribution $\pi$ satisfies the generalized $\rmT_p$ inequality for some increasing function $\rmJ : \RR^+ \rightarrow \RR^+$. 
  \end{assumptionp}
  If we take $J(r) \equiv \sqrt{2r/\lambda_{\rmT}}$,for every $r \in \RR^+$, then we retrieve the classical $\rmT_p$ inequality. Thus, we indeed generalize 
  $\rmT_p$ inequality with this assumption.  
  An important example of \Cref{Tp-assumption} is a consequence of  \cite{bolley2005weighted}[Corollary 2.3]. 
  They prove that if 
  \begin{equation*}
    \int_{\RR^d} \exp(\|x - x_0\|^p) \pi(\rmd x) < +\infty, \text{ for some }x_0 \in \RR^d,
  \end{equation*}
  then \Cref{Tp-assumption} is satisfied for $\rmJ(r) =  \lambda_{BV}(r^{\nicefrac{1}{p}} + (r/2)^{\nicefrac{1}{2p}})$, where $\lambda_{BV}$ is a constant that may depend on the dimension (see \Cref{sec:app-Tp-assumption} for details).
  Therefore, in this case the assumption is essentially a condition on the tails of the target distribution $\pi$.
  In particular, for the generalized normal distribution the verification of this bound is straightforward.

  The third assumption is chosen to essentially restrict ourselves for potentials with (at most) polynomial growth.
  Mathematically, it is expressed as follows:


  \begin{assumptionp}{$(\mathrm{poly},\rmQ)$}\label{poly-assumption}
    For some $p>0$, there exists a polynomial  with positive coefficients such that $\mathrm{ord} ({\rmQ}) = p$  and the following inequality is true:
    \begin{equation*}
      \big\|\nabla V(x)\big\| \leq  {\rmQ}( \|x \|).
    \end{equation*}
  \end{assumptionp}
  Using Taylor formula, one may check that $L$-smooth potentials satisfy this condition with $\rmQ(r) = L r + \|\nabla V(0)\|$.
  We want highlight that the constant $p$ is the same for the assumptions \ref{Tp-assumption} and \ref{poly-assumption}. 
  This will allow us to treat the polynomials of order $p$. 
  For detailed examples of distributions that satisfy our set of assumptions please refer to \Cref{sec:app-applications}.

  We conclude our list of assumptions with a bound on the kernel $k(\cdot,\cdot)$. 
  \begin{assumptionp}{$({\rm ker},B)$}\label{A-ker-bound}
      There exists $B>0$ such that    $\|k(x, .)\|_{\mathcal{H}_{0}} \leq B$ and 
      \begin{equation*}
        \left\|\nabla_{x} k(x, .)\right\|_{\mathcal{H}}=\left(\sum_{i=1}^{d}\left\|\partial_{x_{i}} k\left(x, .\right)\right\|_{\mathcal{H}_{0}}^{2}\right)^{\frac{1}{2}} \leq B,
      \end{equation*}
       for all $x \in \RR^d$.
  \end{assumptionp}
  Due to the reproducing property, these conditions are equivalent to $k(x,x)\leq B$ and $\partial^2_{x_i,x_i} k(x,x) \leq B$. 
  Here, the first and the second partial derivatives are operated, respectively, on the first and the second variables of $k(\cdot,\cdot)$. 
  Based on this criterion, one may check that the multiquadratic kernel $k(x,y) := (c^2 + \|x-y\|^2)^{\beta}$ for some $c>0>\beta>-1$ satisfies the \Cref{A-ker-bound}. See appendices \ref{append:repro} for details on the reproductive property.

\section{Main result}\label{sec:main}

  
  In this section we present our main results. We start with a proposition that bounds the difference between the value of the objective at two consecutive iterations. The proof of the proposition can be found in \Cref{sec:propproof}
  \begin{proposition}\label{prop-descent-lemma} Suppose that Assumptions \ref{A-ker-bound} and \ref{A-l0l1}  are satisfied. Let $\alpha>1$ and choose
    \begin{equation}\label{eq:gamma-cond}
      \gamma \leq   \frac{(\alpha-1)\cdot\min\{ 1 , {1}/{L_1} \} }{\alpha B \|g_{\mu_n}\|_{\mathcal{H}}} .
    \end{equation} 
    Then
    \begin{equation}\label{eq:descent-lemma}
      \mathrm{KL}\left(\mu_{n+1} \mid \pi\right)-\mathrm{KL}\left(\mu_{n} \mid \pi\right)
      \leq - \gamma  \left[ 1 - \frac{\gamma}{2} B^2 \big(\alpha^2 + (e-1) A_n \big) \right] I_{\Stein}(\mu_n \mid \pi),
    \end{equation}
    where $A_n := L_{0}+L_{1}\EE_{X \sim \mu_n}[\|\nabla V(X)\|]$.
  \end{proposition} 
  \Cref{prop-descent-lemma} may be seen as a descent lemma on the $\KL$ divergence.  Let us develop the condition \eqref{eq:gamma-cond}. The following lemma provides us with an upper bound on $\|g_{\mu_n}\|_{\cH}$.
  \begin{lemma}\label{lem:A3impliesA3}
    If \Cref{A-ker-bound} is satisfied, then
    \begin{equation*}
      \|g_{\mu_n}\|_{\cH} = I_{\Stein}(\mu_n\mid \pi)^{\frac{1}{2} }
       \leq B\left(\mathbb{E}_{X \sim \mu_n}\big[\|\nabla V(X)\|\big]+1\right)
    \end{equation*}
    for all $n\in \NN$.
  \end{lemma}
  Thus \eqref{eq:gamma-cond} may be replaced by 
  \begin{equation}\label{eq:gamma-cond-new}
    \gamma \leq (\alpha-1)\min\{ 1 , {1}/{L_1} \}\left[\alpha B^2\left(\mathbb{E}_{X \sim \mu_n}\big[\|\nabla V(X)\|\big]+1\right)\right]^{-1}.
  \end{equation}
  On the other hand, we would like the right-hand side of \eqref{eq:descent-lemma} to be negative and it yields the following condition on the step-size:
  \begin{equation}\label{eq:gamma-cond2}
    \begin{aligned} 
      \gamma &< 2B^{-2}\big({\alpha^2 + (e-1) A_n}\big)^{-1} = 
      { 2B^{-2}\big(\alpha^2 + (e-1) \big(L_{0}+L_{1}\EE_{X \sim \mu_n}[\|\nabla V(X)\|]\big)\big)}^{-1}.
    \end{aligned}
  \end{equation}
  Both in \eqref{eq:gamma-cond-new} and \eqref{eq:gamma-cond2} we have dependence on $\EE_{X \sim \mu_n}[\|\nabla V(X)\|]$. 
  The lemma below establishes an upper bound on this quantity.
  \begin{lemma}\label{lem:exp-grad}
    Suppose that the potential function $V$ satisfies the Assumptions \ref{Tp-assumption} and \ref{poly-assumption} for some constants  $p$ and a polynomial $\rmQ$. Then,
    \begin{equation*}
      \EE_{X \sim \mu_n}[\|\nabla V(X)\|] \leq \rmQ \left( \rmJ\big(\KL(\mu_n\mid \pi)\big) + W_p(\pi,\delta_0)\right).
    \end{equation*}
  \end{lemma}
  The important implication of \Cref{lem:exp-grad} is that both upper bounds on the step-size are inversely proportional to
  $\KL(\mu_n\mid \pi)$. Let us now state the main theorem.
  \begin{theorem}[Descent lemma]
    \label{thm:main}
    Let the target distribution $\pi$ and its potential function $V$ satisfy the Assumptions \ref{A-ker-bound}, \ref{A-l0l1},
    \ref{Tp-assumption}, and \ref{poly-assumption}. 
    Define $C_0 := \rmQ \left( \rmJ\big(\KL(\mu_0\mid \pi)\big) + W_p(\pi,\delta_0)\right)$ and  suppose that the step-size $\gamma$ satisfies 
    \begin{equation}\label{eq:gamma-short}
      \gamma < \frac{\alpha - 1}{ \alpha B^2 \big(\alpha^2 + (e - 1) (\max (L_0,L_1,1) + \max (L_1,1) C_0)\big)}.
    \end{equation}
    Then for every $n = 0,1,\ldots$ the following inequality is true: 
    \begin{equation*}\label{eq:thm-main}
      \mathrm{KL}\left(\mu_{n+1} \mid \pi\right)-\mathrm{KL}\left(\mu_{n} \mid \pi\right)
      \leq - \frac{\gamma}{2} I_{\Stein}(\mu_n \mid \pi).
    \end{equation*}
  \end{theorem}
  The proof is postponed to \Cref{sec:app-thm-proof}. 
  
  \begin{corollary}[Convergence]\label{cor:convergence}
    Let the assumptions of \Cref{thm:main} be satisfied.
    Then following statements are true.
    \begin{itemize}
      \item[1.] $I_{\Stein} (\mu_n \mid \pi)$   converges to zero, when $n\rightarrow\infty$. 
      \item[2.] The average Stein-Fisher error of the first $n$ iterates converges to 0, with $O(1/n)$ rate:
      \begin{equation*}
        \frac{1}{n} \sum_{i=0}^{n-1} I_{\Stein}\left(\mu_{i} \mid \pi\right) 
        \leq \frac{2}{n \gamma}{\KL (\mu_0\mid \pi) }.
      \end{equation*}
    \end{itemize}
  \end{corollary}
  This yields in particular that the series with the general term $I_{\Stein}(\mu_n \mid \pi)$ is convergent, that is $I_{\Stein}(\mu_n \mid \pi) \rightarrow 0$.
  When $\pi$ is distant dissipative and the $k$ is the inverse multiquadratic kernel, then convergence in $I_{\Stein}$ yields weak convergence of the sequence $\mu_n$ to the target $\pi$ (see \Cref{sec:app-weak-conv}).
  Proofs of the corollary can be found in \Cref{sec:app-proof-cor-conv}.  
  \Cref{cor:convergence} yields that if  $n > \nicefrac{\KL(\mu_0\mid\pi)}{\gamma\varepsilon}$,
  then $\frac{1}{n} \sum_{i=0}^{n-1} I_{\Stein}\left(\mu_{i} \mid \pi\right)<\varepsilon$. 
  Based on this inequality, the next theorem estimates the complexity of the algorithm for some particular choices of $\rmJ$ and $\mu_0$.
  \begin{theorem}\label{thm:complexity}
    Let assumptions \ref{A-ker-bound}, \ref{A-l0l1}, and \ref{poly-assumption} hold and let $\mu_0=\cN(0,{\rm I}_d)$. 
    Then in order to have $\varepsilon$ average Stein-Fisher error it is sufficient to perform 
    $n$ iterations of the SVGD, where
    \begin{itemize}
      \item $n=\cO \left({\varepsilon}^{-1}{Q(1)^3 \max(L_1,1)\lambda_{BV}^p (pd)^{p+1}}\right)$, 
      if the target $\pi$ satisfies \Cref{Tp-assumption}  with $\rmJ(r) = \lambda_{BV} (r^{\nicefrac{1}{p}} + (r/2)^{\nicefrac{1}{2p}})$;
      \item  $n=\cO \left({\varepsilon}^{-1} { \rmQ(1)^{\frac{p+2}{2}} \max(L_1,1)\lambda_{\rmT}^{-\frac{p}{2}} \left(p d\right)^{\frac{(p+1)(p+2)}{4}} }\right)$,
      if the target $\pi$ satisfies Talagrand's $\rmT_p$ inequality \eqref{eq:Tp} with constant $\lambda_{\rmT}$.
    \end{itemize}
  \end{theorem}
  The proof is postponed to \Cref{sec:app-proof-complexity}. 
  We observe that in both cases we achieve polynomial convergence of the algorithm in terms of the dimension and the precision. 
  As mentioned in \Cref{sec:setting}, in the first setting, \ref{Tp-assumption} essentially becomes a tail condition, which can be easily verified for a wide class of log-polynomial densities. The Talagrand's $\rmT_p$ is harder to verify in the general case. 
  However, unlike $\lambda_{\rmT}$, the constant $\lambda_{BV}$ may be dimension dependent (see \Cref{sec:app-Tp-assumption}). 
  Therefore, when $p \leq 2$ Talagrand's inequality yields better convergence. 
  In particular, when $p = 1$, we recover the complexity $\cO(d^{3/2})$ which coincides with the result from 
  \citep{salim2021complexity}.
  On the other hand, when $p$ is large, the algorithm is likely to perform significantly faster under the generalized $\rmT_p$.
\section{Examples}\label{sec:app-applications}
  In this section, we discuss a few cases that our method can be applied on.
  As mentioned in \Cref{sec:setting}, the family of generalized Gaussian distributions satisfy our set of assumptions:
  \begin{equation*}
    \pi(x)\propto \exp\left(-\frac{\|x - a\|^p}{2\sigma^p}\right).
  \end{equation*}
  This family of distributions has received considerable attention from the engineering community, due to the flexible parametric form of its 
  probability density function in modeling many physical phenomena. 
  For instance values of $p=2.2$ and $p=1.6$ have been found to model the ship transit noise and the sea surface agitation noise respectively \cite{banerjee2013underwater}. 
  Due to its simple explicit form of the density function, many important probabilistic quantities, such as moment, entropy etc. are easy to compute for the generalized Gaussians. It is known that that this class minimizes the entropy under  under a $p$-th absolute moment constraint \citep{cover2006elements}.
  An overview of the analytic properties of this family of distributions can be found in \cite{dytso2018analytical}. 
  As seen previously, when $p\geq 2$, this family of distributions meets the assumptions of this work, so we can apply SVGD on them.

  The second example has its roots in the Bayesian LASSO problem . 
  The LASSO estimator (proposed by\citep{tibshirani1996regression}) can be seen as the $\ell_1$-penalized least square estimator:
  This can be interpreted also as a Bayesian posterior mode (also known as Bayesian LASSO), 
  with Laplacian priors (see \citep{park2008bayesian,fu1998penalized,mallick2014new}).

  From the sampling perspective, the Bayesian LASSO boils down to 
  \begin{equation*}
    \pi(x)\propto \exp(-f(x)-\tau\|x\|_p^p).
  \end{equation*}
  Here $f$ is the likelihood function of the prior. In \citep{park2008bayesian}, the authors consider the simpler case where $f(x)=\|{A x}-{b}\|^{2}$, with 
  ${A}$ is the measurement matrix, while ${b}$ is the vector of labels. 
  In the case when $p > 1$ the classical sampling methods fail to provide guarantees as the potential does not satisfy the smoothness condition. 
  When $p\in(1,2)$, \cite{chatterji-lmc-smoothness} applied the Gaussian smoothing technique 
  on the potential function $V(x)=f(x)+\tau\|x\|_p^p$. 
  This allowed to perform the usual Langevin Monte-Carlo algorithm on the modified potential. 
  However, their method does not deal with the case $p\geq 2$. Moreover, they assume the strong convexity of the potential. 
  Our analysis for the SVGD can deal with the case when $p\geq 2$. In particular, we do not need to put any restriction on the matrix ${A}$.

\section{Conclusion and Future work}\label{sec:conclusion}

In this work, we studied the convergence of the SVGD under a relaxed smoothness condition. 
The latter was first proposed in the context of gradient clipping for the optimization problem.
This relaxed smoothness condition allows to treat distributions with polynomial potentials. 
Main result consists of a descent lemma for the $\KL$  error of  the SVGD algorithm.
The result implies polynomial convergence for the average Stein-Fisher information error for two types of functional inequalities. 
However, at the moment this setting is of purely theoretical interest as the algorithm operates on the space of probability distributions.
Despite numerous applications, the discretized algorithm does not have practical convergence guarantees in the general case. 
This problem remains open for future work.

\bibliographystyle{plainnat}
\bibliography{bibliography.bib}

\newpage
\addtocontents{toc}{\protect\setcounter{tocdepth}{2}}

\appendix
\part*{Appendix}
\tableofcontents
\newpage

\section{Langevin Monte-Carlo from the optimization perspective}\label{sec:app-langevin}

  Langevin algorithms have extensively been studied by statistics and machine learning communities during the last decades. 
  The algorithms rely on an SDE called (vanilla) Langevin dynamics:
  \begin{equation*}
    \label{eq:LangevinDynamic}
    dX_t=-\nabla V(X_t)dt+\sqrt{2}dB_t,
  \end{equation*}
  where $B_t$ is a $d$-dimensional Brownian motion on $\RR^d$. 
  The important property of this SDE is that under technical conditions the target $\pi(\cdot) \propto \exp(-V(\cdot))$ is its unique invariant distribution. 
  Moreover, it is ergodic and $\KL(\mu_t,\pi)$ converges linearly to zero, where $\mu_t$ is the distribution of $X_t$ \citep{bhattacharya1978criteria}.
  The Euler-Maruyama discretization of this SDE over the time axis results in the Langevin Monte-Carlo (LMC) 
  algorithm:
  \begin{equation}\label{eq:LMC}
    \theta_{i+1} = \theta - \gamma \nabla V (\theta_i) + \sqrt{2h} \xi_{i+1}, \text{ for } i\in \NN.
  \end{equation}
  Here $\gamma>0$ is the discretization step-size and $\xi_{1},\xi_{2},\ldots$ is a sequence of independent standard Gaussians independent from the initial point $\theta_0$. 
  This algorithm initially was studied by \cite{roberts1996exponential,roberts1998optimal,roberts2002langevin,dwivedi2018log} who proposed to apply Metropolis-Hastings step (MALA). 
  The reason for this adjustment is that the constant step-size $\gamma$ induces a bias and the target distribution $\pi$ is no 
  longer the invariant distribution of the discrete-time process.
  Later, \cite{dalalyan2017theoretical} suggested to remove the Metropolis-Hastings step with a result that essentially controls the bias depending on the step-size $\gamma$. This instigated a new line of research which studied the convergence properties of the LMC in different settings (see \cite{durmus2017nonasymptotic,cheng2018sharp,cheng2018convergence,dalalyan2019user,durmus2019high,vempala2019rapid}).
  
  The LMC algorithm essentially performs a Forward-Flow iteration in the space of distributions (see \citep{wibisono2018sampling,durmus2018analysis}). Indeed, if we define by $\mu_i$ the distribution of the $i$-th iterate $\theta_i$, then
  \eqref{eq:LMC} is equivalent to
  \begin{equation*}
    \begin{aligned}
      &\mu_{i+\nicefrac{1}{2}} = \big( I - \gamma \nabla V(\mu_i)\big){\#}\mu_i\\
      &\mu_{i+1} =  \mathcal{N}(0,{\rm I}_d) \hspace{0.1cm} \textasteriskcentered 
      \hspace{0.1cm} \mu_{i+\nicefrac{1}{2}},
    \end{aligned}
  \end{equation*}
  where the $\#$ is the pushforward measure, while $\textasteriskcentered$ is the convolution. 
  The first equation corresponds to the gradient descent for $\EE_{X\sim\mu}[V(X)]$ where the argument is the measure 
  $\mu$. The second equation can be interpreted as the exact gradient flow of the negative entropy. 
  The combination of this two steps results in a biased algorithm (as noticed in \citep{roberts1996exponential}) since the flow step is not the adjoint of the gradient step.

  \section{Convergence results}\label{sec:mainproofs}

  \subsection{Proofs of the \Cref{prop-descent-lemma}}\label{sec:propproof}

    For the sake of brevity we will omit the index when referring to $g_{\mu_n}$ and will write simply $g$.
    Let us define $\phi_{t} := I - t g$ and  $\rho_t:= {\phi_{t}}_{\#} \mu_n$ for every $t \in[0, \gamma]$.  
    Then, applying Taylor formula to the function 
    \begin{equation}\label{eq:g9y980fd8}
      \varphi(t):=\mathrm{KL}\left(\rho_{t} \mid \pi\right)
    \end{equation} 
    we have the following
    \begin{equation}\label{eq:phi-taylor}
      \varphi(\gamma)=\varphi(0)+\gamma \varphi^{\prime}(0)+\int_{0}^{\gamma}(\gamma-t) \varphi^{\prime \prime}(t) {\rmd} t.
    \end{equation}
    By the definition of the SVGD iteration, we have that $\varphi(0)={\KL}\left(\mu_{n} \mid \pi\right)$ and $\varphi(\gamma)={\KL}\left(\mu_{n+1} \mid \pi\right)$. Let us now compute the term of \eqref{eq:phi-taylor} corresponding to the first order derivative.
    \begin{lemma}\label{lem:diff}
      Suppose that Assumption \ref{A-ker-bound} holds. Then, for any $x \in \RR^d$ and  $h \in \mathcal{H}$, 
      \begin{equation*} 
         \|J h(x)\|_{{\rm HS}} \leq B\|h\|_{\mathcal{H}}.
      \end{equation*}
    \end{lemma}
    The proof of the lemma can be found in \Cref{proof:lem:diff}. 
    Applying the lemma to the function $g$, we obtain the following:
    \begin{equation}\label{eq:ineq-op-HS}
      \|t J g(x)\|_{op} \leq \|t J g(x)\|_{{\rm HS}} \leq t B\|g\|_{\mathcal{H}} < 1.
    \end{equation}
    The latter inequality is due to the condition on the step-size $\gamma$. The bound \eqref{eq:ineq-op-HS} implies
    that $\phi_{t}$ is a diffeomorphism. 
    Therefore, $\rho_{t}$ admits a density given by the change of variables formula:
    \begin{equation*}
      \rho_{t}(x)=\left|J \phi_{t}\left(\phi_{t}^{-1}(x)\right)\right|^{-1} \mu_{n}\left(\phi_{t}^{-1}(x)\right).
    \end{equation*}
    Changing the variable of integration and applying the transfer lemma we get the following formula for $\varphi(t)$:
    \begin{equation*}
      \begin{aligned}
        \varphi(t) &=\int \log \left(\frac{\rho_{t}(y)}{\pi(y)}\right) \rho_{t}(\rmd y) \\
        &=\int \log \left(\frac{\mu_{n}(x)\left|J \phi_{t}(x)\right|^{-1}}{\pi\left(\phi_{t}(x)\right)}\right) \mu_{n}(\rmd x)\\
        &=\int  \left[ \log \left(\mu_{n}(x)\right) + \log \left(\left|J \phi_{t}(x)\right|^{-1}\right) - \log \left({\pi\left(\phi_{t}(x)\right)}\right)
        \right]  \mu_{n}(\rmd x).
      \end{aligned}
    \end{equation*}
    Let us then compute the time derivative of $\varphi(t)$. Taking the derivative inside and applying Jacobi's formula for matrix determinant differentiation we obtain the following equality: 
    \begin{equation*}
      \varphi^{\prime}(t)=-\int \operatorname{tr}\left((J \phi_{t}(x))^{-1} \frac{d J \phi_{t}(x)}{d t}\right) \mu_{n}(\rmd x)
      - \int\left\langle\nabla \log \pi\left(\phi_{t}(x)\right), \frac{d \phi_{t}(x)}{d t}\right\rangle \mu_{n}(\rmd x).
    \end{equation*}
    By definition, $d\phi_t / dt = g$. Therefore, we can use the explicit expression of $\phi_{t}$ to write:
    \begin{equation}\label{eq:phi-prime}
      \varphi^{\prime}(t)=\int \operatorname{tr}\left((J \phi_{t}(x))^{-1} J g(x)\right) \mu_{n}(\rmd x)
      + \int \left \langle \nabla V\left(\phi_{t}(x)\right), g(x)\right\rangle \mu_{n}(\rmd x) .
    \end{equation}

    
    The Jacobian at time $t=0$ is simply equal to the identity matrix since $\phi_{0}={\rm I}_d$.
   It follows that $\operatorname{tr}\left((J \phi_{0}(x))^{-1} J g(x)\right)=$ $\operatorname{tr}(J g(x))=\operatorname{div}(g)(x)$ by the definition of the divergence operator. Using integration by parts:
    \begin{equation*}
      \begin{aligned}
        \varphi^{\prime}(0) &=-\int[-\operatorname{div}(g)(x)-\langle\nabla \log \pi(x), g(x)\rangle] \mu_{n}(\rmd x)\\
        &=-\int\left\langle\nabla \log \left(\frac{\mu_{n}}{\pi}\right)(x), g(x)\right\rangle \mu_{n}(\rmd x).
      \end{aligned}
    \end{equation*}
     Based on the alternative definition of $g_{\mu}$ (see \Cref{rem:alt-def-gmu}) and the reproducing property, we have 
    \begin{equation*}
      \begin{aligned}
        \int\Big\langle\nabla\log\left(\frac{\mu_n}{\pi}\right)(x)&,g(x)\Big\rangle\mu_n(\rmd x)\\
        &=\iint k(x,y)\left\langle\nabla\log\left(\frac{\mu_n}{\pi}\right)(x),
        \nabla\log\left(\frac{\mu_n}{\pi}\right)(y)\right\rangle \mu_n(\rmd x)\mu_n(\rmd y)\\
        &=\iint \langle k(x,\cdot),k(y,\cdot)\rangle_{\cH_0}\left\langle \nabla\log\left(\frac{\mu_n}{\pi}\right)(x),
        \nabla\log\left(\frac{\mu_n}{\pi}\right)(y)\right\rangle \mu_n(\rmd x)\mu_n(\rmd y)\\
        &=\left\langle \int\nabla\log\left(\frac{\mu_n}{\pi}\right)(x)k(x,\cdot)\mu_n(\rmd x),
        \int\nabla\log\left(\frac{\mu_n}{\pi}\right)(y)k(y,\cdot)\mu_n(\rmd y)\right\rangle_{\cH}\\
        &=\norm{g}_{\cH}^2. 
      \end{aligned}
    \end{equation*}
    Therefore,
    \begin{equation}\label{eq:phi'(0)}
      \varphi^{\prime}(0) = - \norm{g}_{\cH}^2.
    \end{equation}

    Next, we calculate the term of \eqref{eq:phi-taylor} that contains the second derivative using \eqref{eq:phi-prime}.
    First, 
    \begin{equation*}
        \begin{aligned}
          \frac{d}{dt}\operatorname{tr}\left((J \phi_{t}(x))^{-1} J g(x)\right) 
          &=\operatorname{tr}\left(\frac{d}{dt}(J \phi_{t}(x))^{-1} J g(x)\right) \\
          &=-\operatorname{tr}\left((J \phi_{t}(x))^{-2} \frac{d}{dt}J\phi_t(x)J g(x)\right) \\
          &=\operatorname{tr}\left(J \phi_{t}(x)^{-2} J g(x)J g(x)\right) \\
          &=\operatorname{tr}\left(\left(J \phi_{t}(x)^{-1} J g(x)\right)^2\right).
        \end{aligned}
    \end{equation*}
    From the definition of the function $\phi_t$, we know that $J \phi_{t}(x) = ({\rm I}_d + t Jg)(x)$. 
    Thus, $(J \phi_{t})^{-1}$ and $Jg$ commute. The latter yields
    $\operatorname{tr}\left(\left(J g(x)\left(J \phi_{t}(x)\right)^{-1}\right)^{2}\right)
    =\left\|J g(x)\left(J \phi_{t}(x)\right)^{-1}\right\|_{\rm H S}^{2}$. 
    On the other hand, applying the chain rule on second term of \eqref{eq:phi-prime} yields
    \begin{equation*}
      \partial_t \big(\left\langle \nabla V\left(\phi_{t}(x)\right), g(x)\right\rangle\big)
      =   \left\langle g(x), {\nabla^2 V}\left(\phi_{t}(x)\right) g(x)\right\rangle.
    \end{equation*}
    Summing up, we have the following:

        
    \begin{equation*}
      \varphi^{\prime \prime}(t)
      =\underbrace{\int \left\|J g(x)\left(J \phi_{t}(x)\right)^{-1}\right\|_{\rm H S}^{2} \mu_{n}(\rmd x)}_
      {:=\psi_1(t)}
      +\underbrace{\int\left\langle g(x), {\nabla^2 V}\left(\phi_{t}(x)\right) g(x)\right\rangle\mu_{n}(\rmd x)}_{:=\psi_2(t)} .
    \end{equation*}
    First, we bound $\psi_{1}(t)$. Cauchy-Schwarz implies that
    \begin{equation*}
      \left\|J g(x)\left(J \phi_{t}(x)\right)^{-1}\right\|_{\rm H S}^{2} 
      \leq \|J g(x)\|_{\rm H S}^{2}\left\|\left(J \phi_{t}(x)\right)^{-1}\right\|_{o p}^{2}.
    \end{equation*}
    From \Cref{lem:diff}, we have  $\|J g(x)\|_{\rm H S} \leq B\|g\|_{\cH}$. 
    To bound the second term, let us recall that $\phi_t = I - tg$ and that $t \leq \gamma$. Thus, the following bound is true:
    \begin{equation*}
      \begin{aligned}
        \left\|\left(J \phi_{t}(x)\right)^{-1}\right\|_{o p} 
        & = \left\|\big( ( {\rm I}_d - t J g)(x)\big)^{-1}\right\|_{o p} \leq \sum_{i=0}^{\infty}\|t J g(x)\|_{o p}^{i} \leq \sum_{i=0}^{\infty}\|\gamma J g(x)\|_{\rm H S}^{i}.
      \end{aligned}
    \end{equation*}
    Recalling \eqref{eq:gamma-cond} and combining it with \Cref{lem:diff} we obtain 
    \begin{equation*}
      \begin{aligned}
        \left\|\left(J \phi_{t}(x)\right)^{-1}\right\|_{o p} \leq \sum_{i=0}^{\infty}(\gamma B \|g\|_{\cH})^{i} \leq \sum_{i=0}^{\infty} \left(\frac{\alpha-1}{\alpha}\right)^i  = \alpha.
      \end{aligned}
    \end{equation*}
    Summing up, we have that
    \begin{equation*}\label{eq:psi1-bound}
      \psi_1(t) \leq \alpha^2 B^2\|g\|_{\cH}^2.
    \end{equation*}
    Next, we bound $\psi_{2}(t)$. By definition,
    \begin{align*}
      \psi_2(t) &= \EE_{X \sim \mu_n} \left[\left\langle g(X), {\nabla^2 V}\left(\phi_{t}(X)\right) g(X)\right\rangle \right] \leq \EE_{X \sim \mu_n} \left[ \|{\nabla^2 V}\left(\phi_t(X)\right)\|_{op} \|g(X)\|_2^2 \right].
    \end{align*}
    Let us bound the norm of $g(x)$. The reproduction property of the RKHS yields the following:  
    \begin{equation}\label{eq:norm2-H}
      \|g(x)\|_2^{2}=\sum_{i=1}^{d}\left\langle k(x, .), g_{i}\right\rangle_{\mathcal{H}_{0}}^{2} \leq\|k(x, .)\|_{\mathcal{H}_{0}}^{2}\|g\|_{\mathcal{H}}^{2} \leq B^{2}\|g\|_{\mathcal{H}}^{2}.
    \end{equation}
    Therefore,
    \begin{equation}\label{eq:psi2-initial-bound}
      \psi_2(t) \leq  B^{2}\|g\|_{\mathcal{H}}^{2}\EE_{X \sim \mu_n} \left[ \|{\nabla^2 V}\left(\phi_t(X)\right)\|_{op} \right].
    \end{equation}
    Let us bound $\EE_{X \sim \mu_n} \left[ \|{\nabla^2 V}\left(\phi_t(X)\right)\|_{op} \right]$. \Cref{A-l0l1} implies the following inequality:

    \begin{equation*}
        \left\| \nabla^2 V (\phi_t(x))\right\|_{op} \leq L_{0}+L_{1}\|\nabla V(\phi_t(x))\|,
    \end{equation*}
    for every $x \in \RR^d$.
    To bound the term $\|\nabla V(\phi_t(x))\|$, we introduce the following lemma.

    \begin{lemma}\label{lem:grad-v}
      Let $V$ be an $\left(L_{0}, L_{1}\right)$-smooth function and $\Delta>0$ be a constant. 
      For any $x,x^{+} \in \RR^d$ such that $\left\|x^{+}-x\right\| \leq \Delta $, we have 
      \begin{equation*}
        \left\|\nabla V\left(x^{+}\right)\right\| \leq 
        \frac{L_0 }{ L_1 } \left(\exp(\Delta L_1) - 1\right)  +  \|\nabla V(x)\|\exp(\Delta L_1).
      \end{equation*}
    \end{lemma} 
    We will apply \Cref{lem:grad-v} to $\phi_t(x)$ and $\phi_0(x)$. By definition, $\phi_t(x) - \phi_0(x) = t g(x)$ and
    according to inequality \eqref{eq:norm2-H},
    \begin{equation*}
      \|\phi_t(x) - \phi_0(x)\|_2\leq t B\| g\|_{\cH}.
    \end{equation*}
    Thus, using \Cref{lem:grad-v} for $x = \phi_0(x) $, $x^+ = \phi_t(x)$ and $\Delta = tB\|g\|_{\cH}$, we obtain the following:
    \begin{equation}\label{eq:bound-Hess-phi}
      \begin{aligned}
        \left\| \nabla^2 V (\phi_t(x))\right\|_{op} 
        &\leq  L_0 + L_1 \left(\frac{L_0}{L_1}\big(\exp(t B \|g\|_{\cH} L_1) - 1\big) 
        + \left\|\nabla V\left(\phi_0(x) \right)\right\| \exp(t B \|g\|_{\cH}L_1)\right) \\
        &= \big(L_{0} + L_1\left\|\nabla V\left(x\right)\right\|\big)\exp (t B L_1  \|g\|_{\cH} ),
      \end{aligned}
    \end{equation}    
    Combining \eqref{eq:psi2-initial-bound} and \eqref{eq:bound-Hess-phi} we obtain
    \begin{equation*}\label{eq:psi2-final-bound}
      \psi_2(t)\leq   B^2 \|g\|_{\cH}^2 \left(L_{0} + L_1 \EE_{X \sim \mu_n}\big[\left\|\nabla V (X) \right\|\big]\right) 
      \exp (t B L_1 \|g\|_{\cH} ).
    \end{equation*}
    Summing up, the bounds on $\psi_1$ and $\psi_2$ yield the following inequality:
    \begin{equation*}
      \varphi^{\prime \prime}(t) \leq B^2 \|g\|_{\cH}^2 \Big[ \alpha^2 
      + \Big(L_{0} + L_1 \EE_{X \sim \mu_n}\big[\left\|\nabla V (X) \right\|\big]\Big) 
      \exp (t B L_1\|g\|_{\cH} )\Big].
    \end{equation*}
    Recall that by definition $A_n = \big(L_{0} + L_1 \EE_{X \sim \mu_n}\big[\left\|\nabla V (X) \right\|\big]\big).$ 
    Inserting the previous inequality along with \eqref{eq:phi'(0)} to \eqref{eq:phi-taylor}, we get the following bound:
    \begin{equation*}
      \begin{aligned}
        \varphi(\gamma) - \varphi(0) &\leq - \gamma \|g\|_{\cH}^2
        + \Big[ \frac{1}{2}\gamma^2\alpha^2  B^2 \|g\|_{\cH}^2  + A_n  \big( \exp (\gamma B L_1 \|g\|_{\cH}) - \gamma B  L_1\|g\|_{\cH} - 1 \big) 
        \Big].
      \end{aligned}
    \end{equation*}
    One can check that $\exp(t) - t - 1 \leq (e-1)t^2/2$, when $t \in [0,1]$. Since $\gamma B L_1\|g\|_{\cH} < 1$, we deduce
    \begin{equation*}\label{eq:final-prop}
      \varphi(\gamma) - \varphi(0) \leq  \Big[ - \gamma 
        + \frac{1}{2} (\gamma^2 \alpha^2 B^2 + (e-1) A_n  \gamma^2 B^2)\Big] \|g\|_{\cH}^2 .
    \end{equation*}
    Finally, by the definition of Stein's information $\|g\|_{\cH}^2 = I_{\Stein}(\mu_n \mid \pi)$.
    This concludes the proof.
    \begin{remark}
      In fact, the derivation of \eqref{eq:phi'(0)} contains the proof of the fact that $\|g\|_{\cH} = \sqrt{I_{\Stein}(\mu_n\mid\pi)}$.
      Indeed, by the definition of the SVGD (\Cref{sec:SVGD}) the direction function $g$ is chosen to minimize the descent the most.
      On the other hand, the Stein-Fisher information is the maximum of the Stein's linear operator which coincides with our objective.
      Thus, the absolute value of the objective function at $g$, that is $|\phi'(0)|$ 
      equals $\sqrt{I_{\Stein}(\mu_n\mid\pi)}$. 
    \end{remark}
  
  \subsection{Proof of \Cref{thm:main}}\label{sec:app-thm-proof}
  \begin{proof}
    Let us first prove that the sequence $\KL (\mu_n \mid \pi)$ is monotonically decreasing. We will use the method of 
    mathematical induction. First let us notice that \eqref{eq:gamma-short} implies the following system of inequalities for the step-size:
   \begin{equation} \label{eq:gamma-final-cond}
      \begin{cases}
        \gamma &\leq 
        {(\alpha-1)}\min\{ 1 , {1}/{L_1} \}\left[{\alpha B^2\left( C_0 + 1 \right)}\right]^{-1};\\
        \gamma &\leq 
          B^{-2}\Big[\alpha^2 + (e-1) \Big(L_{0}+ L_1 C_0\Big) \Big]^{-1}.
      \end{cases}
    \end{equation}
    For $n=0$ the first equation of the system \eqref{eq:gamma-final-cond} combined with \Cref{lem:exp-grad} 
    implies that the step-size condition is satisfied. Thus, according to \Cref{prop-descent-lemma},
    \begin{equation*}
      \mathrm{KL}\left(\mu_{1} \mid \pi\right)-\mathrm{KL}\left(\mu_{0} \mid \pi\right)
      \leq - \gamma  \left[ 1 - \frac{\gamma}{2} B^2(\alpha^2  + (e-1) (L_0 + L_1 C_0))\right] I_{\Stein}(\mu_0 \mid \pi).
    \end{equation*}
    On the other hand, from the second inequality of \eqref{eq:gamma-final-cond} we get that 
    \begin{equation*} 
      1 - \frac{\gamma}{2} B^2(\alpha^2  + (e-1) (L_0 + L_1 C_0))\geq \frac{1}{2}.
    \end{equation*}
    Therefore, $\mathrm{KL}\left(\mu_{1} \mid \pi\right) \leq \mathrm{KL}\left(\mu_{0} \mid \pi\right)$. 
    Let us define by $C_n:= \rmQ \left( \rmJ\big(\KL(\mu_n\mid \pi)\big) + W_p(\pi,\delta_0)\right)$.
    Since $\rmQ$ and $\rmJ$ are monotonically increasing for positive arguments, we obtain $C_1 \leq C_0$. Therefore,
    \begin{equation*}
      \begin{cases}
        \gamma &\leq 
        {(\alpha-1)}\min\{1,{1}/{L_1}\}\left[{\alpha B^2\left( C_1 + 1 \right)}\right]^{-1};\\
        \gamma &\leq 
          B^{-2}\Big[\alpha^2 + (e-1) \big(L_{0}+ L_1 C_1\big)\Big]^{-1}.
      \end{cases}
    \end{equation*}
    We retrieve \eqref{eq:gamma-final-cond} where the term $C_0$ is replaced by $C_1$.
    Thus, we can repeat the previous arguments for $\mu_1$ and $\mu_2$.  Similarly, we can iterate till $n$. 
    Therefore we obtain the following descent bound:
    \begin{equation*}
      \begin{aligned}
        \mathrm{KL}\left(\mu_{n+1} \mid \pi\right)-\mathrm{KL}\left(\mu_{n} \mid \pi\right)
        &\leq - \gamma  \left[ 1 - \frac{\gamma}{2} B^2 (\alpha^2   + (e-1) (L_0 + L_1 C_n))\right] I_{\Stein}(\mu_n \mid \pi).
      \end{aligned}
    \end{equation*}
    This also yields that the sequence $\KL(\mu_n)$ is decreasing. Since $\rmJ$ and $\rmP$ are monotonically increasing for positive arguments, 
    $C_n$ is also decreasing. Thus,
    \begin{equation*}
      \begin{aligned}
        \mathrm{KL}\left(\mu_{n+1} \mid \pi\right)-\mathrm{KL}\left(\mu_{n} \mid \pi\right)
        &\leq - \gamma  \left[ 1 - \frac{\gamma}{2} B^2 (\alpha^2   + (e-1) (L_0 + L_1 C_0))\right] I_{\Stein}(\mu_n \mid \pi)\\
        &\leq -\frac{\gamma}{2}I_{\Stein}(\mu_n \mid \pi).
      \end{aligned}
    \end{equation*}
    The last inequality is due to the second condition on the step-size $\gamma$. 
  \end{proof}

   \subsection{Proof of \Cref{cor:convergence}}\label{sec:app-proof-cor-conv}
    Let us start with the first statement. 
    Summing the descent bounds of \Cref{thm:main} for $i=0,1,\ldots,n-1$, we obtain the following:
    \begin{equation*}
      \mathrm{KL}\left(\mu_{n} \mid \pi\right)-\mathrm{KL}\left(\mu_{0} \mid \pi\right)
      \leq - \frac{\gamma}{2}  \sum_{i=0}^{n-1}I_{\Stein}(\mu_i \mid \pi).
    \end{equation*}
    Rearranging the terms we get
    \begin{equation*}
      \begin{aligned}
        \sum_{i=0}^{n-1}I_{\Stein}(\mu_i \mid \pi)
         &\leq \frac{2}{\gamma}(\mathrm{KL}\left(\mu_{0} \mid \pi\right)-\mathrm{KL}\left(\mu_{n} \mid \pi\right))\\
        &\leq \frac{2}{\gamma}\mathrm{KL}\left(\mu_{0} \mid \pi\right).
      \end{aligned}
    \end{equation*}
    This means that the series on the left-hand side is convergent and thus, its general term converges to zero. Thus the first point is proved. Dividing both sides of the previous inequality on $n$, we deduce the second statement.

\section{Complexity results}

\subsection{Proof of \Cref{thm:complexity}}\label{sec:app-proof-complexity}
\begin{proof}
  From \Cref{cor:convergence}, we know that if \eqref{eq:gamma-short} is satisfied, then
  \begin{equation*}
    \frac{1}{n} \sum_{k=0}^{n-1} I_{\Stein}\left(\mu_{k} \mid \pi\right) \leq\frac{2\KL(\mu_0\mid\pi)}{n\gamma}.
  \end{equation*}
  Let us bound the initial $\KL$ error.
\begin{lemma}\label{lem:upperboundmu0}
  Let \Cref{poly-assumption} hold with some polynomial $\rmQ$ and $\mu_0=\cN(0,{\rm I}_d)$. 
  We then have
  \begin{equation*}
    \KL(\mu_0\mid\pi)\leq \frac{d}{2}\log\left(\frac{1}{2{\Pi} e}\right)+V(0)+{\rmQ(1) d}\sqrt{\frac{2}{\Pi}}
    +\frac{\rmQ(1)(2d)^{\frac{p+1}{2}}\Gamma\big(\frac{p+2}{2}\big)}{\sqrt{\Pi}(p+1)},
  \end{equation*}
  where $\Pi$ is the area of the circle of radius $1$.
\end{lemma}
  Since $p\geq 1$, \Cref{lem:upperboundmu0} implies that 
  \begin{equation}\label{eq:kl-mu0complexity}
    \begin{aligned}
      \KL(\mu_0\mid\pi)
      &= \cO \left(\frac{\rmQ(1)(2d)^{\frac{p+1}{2}}\Gamma\big(\frac{p+2}{2}\big)}{\sqrt{\Pi}(p+1)}\right)\\
      &= \cO \left(\rmQ(1) \left({p d}\right)^{\frac{p+1}{2}} \right).
    \end{aligned}
  \end{equation}
  The second equality is due to Stirling's formula: $\Gamma(r+1) = \cO(\sqrt{2\Pi r} (r/e)^r)$, for every $r > 0$. 
  In order to estimate the order of the step-size, let us compute the order of $C_0$.
  By definition,
  \begin{equation*}
    \begin{aligned}
      C_0 &= \rmQ \left( \rmJ\big(\KL(\mu_0\mid \pi)\big) + W_p(\pi,\delta_0)\right) \\
      &\leq \rmQ \left( \rmJ\big(\KL(\mu_0\mid \pi)\big) + W_p(\pi,\mu_0) + W_p(\mu_0,\delta_0)\right)\\
      &\leq \rmQ \left( 2\rmJ\big(\KL(\mu_0\mid \pi)\big)  + W_p(\mu_0,\delta_0)\right).
    \end{aligned}
  \end{equation*}
  Here we applied the triangle inequality of $W_p$ and the \Cref{Tp-assumption}. 

  Let us now consider the first point of the theorem.
  Using $Q(r) \leq Q(1)(r^p + 1)$  and $\rmJ(r) = \lambda_{BV}(r^{\nicefrac{1}{p}} + (r/2)^{\nicefrac{1}{2p}} )=  \cO(\lambda_{BV} r^{\nicefrac{1}{p}})$, we obtain
  \begin{equation}\label{eq:c_0-complexity}
    \begin{aligned}
      C_0 
      &\leq \rmQ(1) \cdot \left(\left( 2\rmJ\big(\KL(\mu_0\mid \pi)\big)  + W_p(\mu_0,\delta_0)\right)^p + 1\right)\\
      & = \cO\Big(\rmQ(1) \cdot \left(\left( 2\lambda_{BV} \KL\big( \mu_0\mid \pi\big)^{\nicefrac{1}{p}}  + W_p(\mu_0,\delta_0)\right)^p + 1\right)\Big)\\
      & = \cO\Big(\rmQ(1) \cdot \left( (4\lambda_{BV})^p \KL\big( \mu_0\mid \pi\big)  + 2^{p}W_p^p(\mu_0,\delta_0) + 1\right)\Big).
      \end{aligned}
  \end{equation}
  Applying \eqref{eq:kl-mu0complexity} we get
   \begin{equation}
    \begin{aligned}
      C_0 
      & = \cO\left(\rmQ(1) \cdot \left( (4\lambda_{BV})^p   \cO \left(\rmQ(1) \left({pd}\right)^{\frac{p+1}{2}} \right) 
      + 2^p W_p^p(\mu_0,\delta_0)\right)\right)\\
      & = \cO\left(\rmQ(1)^2  (4\lambda_{BV})^p  \left(p d\right)^{\frac{p+1}{2}}
       + \rmQ(1) 2^p W_p^p(\mu_0,\delta_0)\right).
    \end{aligned}
  \end{equation}
  The following lemma is to bound the $(p+1)$-th norm of the standard multivariate Gaussian.
  \begin{lemma}\label{lem:gaussian-moment}
     Let $\mu_0$ be the standard multivariate Gaussian defined on $\RR^d$. Then, for every integer $m \geq 2$
      \begin{equation*}
        \int_{\RR^d} \|{x}\|^{m} \mu_0(\rmd x) 
        \leq  \frac{(2d)^{\frac{m}{2}}}{\sqrt{ \Pi}}  \Gamma\left(\frac{m+1}{2}\right).
      \end{equation*}
  \end{lemma}
  
  If we write down the definition of the Wasserstein distance, then \Cref{lem:gaussian-moment} for $m = p$ yields
  \begin{equation}\label{eq:w_p-mu0}
    W_p^p(\mu_0,\delta_0) = \int_{\RR^d} \|x\|^p \mu_0(\rmd x)
    \leq \frac{(2d)^{\frac{p}{2}} \Gamma\left(\frac{p+1}{2}\right)}{\sqrt{\Pi}} 
    = \cO\left(\left({pd}\right)^{\frac{p}{2}} \right).
  \end{equation}
  This implies
  \begin{equation*}
    \begin{aligned}
    C_0 
    &= \cO\left(\rmQ(1)^2  (4\lambda_{BV})^p \left({pd}\right)^{\frac{p+1}{2}}
     + \rmQ(1) 2^p   \left({pd}\right)^{\frac{p}{2}} \right)\\
    &= \cO\left(\rmQ(1)^2  (4\lambda_{BV})^p \left({pd}\right)^{\frac{p+1}{2}} \right).
    \end{aligned}
  \end{equation*}
  Now let us look back at \eqref{eq:gamma-short}. This inequality yields that \Cref{cor:convergence} is true for
  \begin{equation*}
    \begin{aligned}
      \gamma 
      &= \cO\left(\left\{ B^2 (\alpha^2 + (e-1)(L_0 + L_1 \rmQ(1)^2  (4\lambda_{BV})^p \left({pd}\right)^{\frac{p+1}{2}} ))\right\}^{-1}\right)\\
      &= \cO\left(\left\{ B^2  \max(L_1,1) \rmQ(1)^2  \lambda_{BV}^p \left({pd}\right)^{\frac{p+1}{2}} \right\}^{-1}\right).
    \end{aligned}
  \end{equation*}
  This implies the following equality and concludes the proof of the first point:
  \begin{equation*}
    \begin{aligned}
      n = \frac{\KL(\mu_0\mid \pi)}{2\gamma\varepsilon} &= \cO\left(\frac{1}{\varepsilon}\left[{Q(1)\left({pd}\right)^{\frac{p+1}{2}}}\right]
      \left[
      { B^2 (\alpha^2 + (e-1)(L_0 + L_1 \rmQ(1)^2  (4\lambda_{BV})^p \left({pd}\right)^{\frac{p+1}{2}} ))}\right]\right)\\
      &= \cO \left( B^2  \max(L_1,1) \frac{Q(1)^3 \lambda_{BV}^p (pd)^{p+1}}{\varepsilon}\right).
    \end{aligned}
  \end{equation*}
  For the second point, we will do the same analysis.
  In the case of $\rmT_p$ inequality, the target satisfies \Cref{Tp-assumption} with $\rmJ(r) = \sqrt{2r/\lambda_{\rmT}}$, for some constant $\lambda_{\rmT} > 0$.
  Then, similar to \eqref{eq:c_0-complexity}, we obtain the following complexity for $C_0$:
  \begin{equation*}
    \begin{aligned}
      C_0 
      &\leq \rmQ(1) \cdot \left(\left( 2\rmJ\big(\KL(\mu_0\mid \pi)\big)  + W_p(\mu_0,\delta_0)\right)^p + 1\right)\\
      & = \rmQ(1) \cdot \left(\left( 2\sqrt{2 \KL\big( \mu_0\mid \pi\big)/\lambda_{\rmT}}  + W_p(\mu_0,\delta_0)\right)^p + 1\right)\\\
      & \leq \rmQ(1) \cdot \left( (32/\lambda_{\rmT})^{\frac{p}{2}} \KL\big( \mu_0\mid \pi\big)^{\frac{p}{2}}  + 2^{p}W_p^p(\mu_0,\delta_0) + 1\right).
    \end{aligned}
  \end{equation*}
  Applying \eqref{eq:kl-mu0complexity}
  \begin{equation*}
    \begin{aligned}  
      C_0 
      & = \cO\left(\rmQ(1) \cdot \left( (32/\lambda_{\rmT})^{\frac{p}{2}}   \cO \left(\rmQ(1)^{\frac{p}{2}} \left({pd}\right)^{\frac{p(p+1)}{4}} \right) 
      + 2^p W_p^p(\mu_0,\delta_0)\right)\right)\\
      & = \cO\left(\rmQ(1)^{\frac{p+2}{2}} (32/\lambda_{\rmT})^{\frac{p}{2}} \left(p d\right)^{\frac{p(p+1)}{4}}
       + \rmQ(1) 2^p W_p^p(\mu_0,\delta_0)\right).
    \end{aligned}
  \end{equation*}
  Applying \eqref{eq:w_p-mu0}, we conclude 
  \begin{equation*}
    C_0 = \cO\left(\rmQ(1)^{\frac{p+2}{2}} \lambda_{\rmT}^{-\frac{p}{2}} \left(p d\right)^{\frac{p(p+1)}{4}}\right).
  \end{equation*}
  Thus, in order to have $\varepsilon$ average $I_{\Stein}$ error it is  sufficient to perform $n$ iterations, where
  \begin{equation*}
    \begin{aligned}
      n = \frac{\KL(\mu_0\mid \pi)}{\gamma \varepsilon} 
      &= \cO\left(\frac{\left({p d}\right)^{\frac{p+1}{2}}}{\varepsilon}\cdot B^2 \big(\alpha^2 + (e - 1) (\max (L_0,L_1,1) + \max (L_1,1) C_0)\big) \right)\\
      &= \cO\left(B^2 \max (L_1,1)\frac{ \rmQ(1)^{\frac{p+2}{2}} \lambda_{\rmT}^{-\frac{p}{2}} \left(p d\right)^{\frac{(p+1)(p+2)}{4}} }{\varepsilon}  \right).
    \end{aligned}
  \end{equation*}
  This concludes the proof.

\end{proof}

\section{Proofs of the lemmas}

\subsection{Proof of \Cref{lem:A3impliesA3}}
  \begin{proof}
  Denote $\Phi(x):=k(x,\cdot)\in \mathcal{H}$. Then by the definition of the Stein discrepancy we have the following: 
  \begin{equation*}
    \begin{aligned}
      I_{\Stein}(\mu_n\mid \pi)^{\frac{1}{2}} 
      &=\left\|\mathbb{E}_{X \sim \mu_n} \big[ (\nabla V(X) \Phi(X)-\nabla \Phi(X)) \big]\right\|_{\mathcal{H}}\\ 
      & \leq \mathbb{E}_{X \sim \mu_n}\big[\|\nabla V(X) \Phi(X)-\nabla \Phi(X)\big]\|_{\mathcal{H}}.
     \end{aligned}
  \end{equation*} 
  Applying the triangle inequality and Cauchy-Schwartz inequality, we obtain 
  \begin{equation*}
    \begin{aligned}
      I_{\Stein}(\mu_n\mid \pi)^{\frac{1}{2}} 
      & \leq \mathbb{E}_{X \sim \mu_n}\big[\|\nabla V(X) \Phi(X)\|_{\mathcal{H}}\big]
      + \mathbb{E}_{X \sim \mu_n}\big[\|\nabla \Phi(X)\|_{\mathcal{H}}\big] \\ 
      &=\mathbb{E}_{X \sim \mu_n}\big[\|\nabla V(X)\|\|\Phi(X)\|_{\mathcal{H}}\big] + \mathbb{E}_{X \sim \mu_n}
      \big[\|\nabla \Phi(X)\|_{\mathcal{H}} \big] \\
       & \leq B\left(\mathbb{E}_{X \sim \mu_n}\big[\|\nabla V(X)\|\big]+1\right).
    \end{aligned}
  \end{equation*}
  \end{proof} 
  \subsection{Proof of \Cref{lem:exp-grad}}
    Let the polynomial $\rmQ$ have the following explicit form:
    \begin{equation}\label{eq:q-explicit}
      \rmQ (r) = \sum_{i=0}^{m} a_i r^{p_i}
    \end{equation} 
    for every $r \in \RR$. Here $p = p_0 > p_1 > \ldots > p_m$ and $a_i > 0$. 
    Then \Cref{poly-assumption} yields
    \begin{equation*}
      \begin{aligned}
        \EE_{X \sim \mu_n}\big[ \|\nabla V(X)\| \big] 
        &\leq \EE_{X \sim \mu_n}\big[  \rmP(\|X\|)\big] \\
        &= \EE_{X \sim \mu_n}\left[   \sum_{i=0}^{m} a_i \|X\|^{p_i}\right]\\
        &=  \sum_{i=0}^{m}   a_i W_{p_i}\big(\mu_n,\delta_0 \big)^{p_i},
      \end{aligned}
    \end{equation*}
    where $\delta_0$ is Dirac measure on $\RR^d$ at point $0$. 
    Using the fact that $W_r$ is monotonically increasing w.r.t. to $r$, we have the following
    \begin{equation*}     
     \begin{aligned}
        \EE_{X \sim \mu_n}\big[ \|\nabla V(X)\| \big] 
        &\leq  \sum_{i=0}^{m}   a_i W_{p}\big(\mu_n,\delta_0 \big)^{p_i}\\
        &\leq \sum_{i=1}^{m} a_i \big(W_p(\mu_n,\pi) + W_p(\pi,\delta_0)\big)^{p_i}\\
        &\leq \sum_{i=1}^{m} a_i \left( \rmJ\big(\KL(\mu_n\mid \pi)\big) + W_p(\pi,\delta_0)\right)^{p_i}.\\
        &\leq \rmQ \left( \rmJ\big(\KL(\mu_n\mid \pi)\big) + W_p(\pi,\delta_0)\right).
      \end{aligned}
    \end{equation*}
    The third inequality is due to \Cref{Tp-assumption}. Recalling \eqref{eq:q-explicit}  we conclude the proof.

  \subsection{Proof of \Cref{lem:diff}}\label{proof:lem:diff}
    The proof is based on the reproducing property and Cauchy-Schwarz inequality in the RKHS space. Indeed,
    \begin{equation*}
      \begin{aligned}
        \|J h(x)\|_{\rm H S}^{2} &= \sum_{i, j=1}^{d}\left|\frac{\partial h_{i}(x)}{\partial x_{j}}\right|^{2}\\
        &=\sum_{i, j=1}^{d}\left\langle\partial_{x_{j}} k(x, .), h_{i}\right\rangle_{\mathcal{H}_{0}} \\
        &\leq \sum_{i, j=1}^{d}\left\|\partial_{x_{j}} k(x, .)\right\|_{\mathcal{H}_{0}}^{2}\left\|h_{i}\right\|_{\mathcal{H}_{0}}^{2} \\
        &=\|\nabla k(x, .)\|_{\mathcal{H}}^{2}\left\|h\right\|_{\mathcal{H}}^{2} \\
        &\leq B^{2}\left\|h\right\|_{\mathcal{H}}^{2}   .
      \end{aligned}
    \end{equation*}
    This concludes the proof.
  
  \subsection{Proof of \Cref{lem:grad-v}}
  \begin{proof} 
    Let us fix $x,x^+ \in \RR^d$ and let $\tau(t)$ be defined as $\tau(t)=t\left(x^{+}-x\right)+x$ with $t \in[0,1]$. Then we have
    \begin{equation*}
      \nabla V(\tau(t))=\int_{0}^{t} \nabla^{2} V(\tau(u))\left(x^{+}-x\right) \mathrm{d} \tau+\nabla V(\tau(0)) .
    \end{equation*}
    Let us bound the norm of $\nabla V(\tau(t))$. Applying triangle and Cauchy-Schwartz inequalities we obtain:
    \begin{equation*}
      \begin{aligned}
        \|\nabla V(\tau(t))\| & \leq \int_{0}^{t}\left\|\nabla^{2} V(\tau(u))\left(x^{+}-x\right)\right\| 
        \mathrm{d}u +\|\nabla V(\tau(0))\| \\
        & \leq\left\|x^{+}-x\right\| \int_{0}^{t}\left\|\nabla^{2} V(\tau(u))\right\| \mathrm{d} u
         + \|\nabla V(x)\|.
      \end{aligned}
    \end{equation*}
    \Cref{A-l0l1} yields
    \begin{equation*}
      \begin{aligned}
         \|\nabla V(\tau(t))\| 
         &\leq  \Delta \int_{0}^{t}\left(L_{0}+L_{1}\|\nabla V(\tau(u))\|\right) \mathrm{d} u+\|\nabla V(x)\|\\
          & =  \Delta L_{0}t +  \|\nabla V(x)\| + \int_{0}^{t} \Delta L_{1}\big\|\nabla V(\tau(u))\big\| \mathrm{d} u.
      \end{aligned}
    \end{equation*}
    Applying Grönwall's integral inequality (see \Cref{lem:gronwall}) to the function $ \|\nabla V(\tau(\cdot))\| $  we obtain
    \begin{equation*}
      \begin{aligned} 
         \|\nabla V(\tau(t))\| 
          & =  \Delta L_{0}t +  \|\nabla V(x)\| 
          + \int_{0}^{t} \Delta L_{1}\big(\Delta L_{0}u  +  \|\nabla V(x)\|\big)\exp\big(\Delta L_{1}(t-u)\big)\mathrm{d} u.
      \end{aligned}
    \end{equation*}
    Inserting $t = 1$, we get the following:
    \begin{equation*}
      \begin{aligned}
         \|\nabla V(x^+)\| 
         & \leq \Delta L_{0} +  \|\nabla V(x)\| 
          + \int_{0}^{1} \Delta L_{1}\big(\Delta L_{0}u  +  \|\nabla V(x)\|\big)\exp\big(\Delta L_{1}(1-u)\big)\mathrm{d} u\\
          & = - \frac{L_0}{L_1} + \left(\frac{L_0 }{ L_1 }  +  \|\nabla V(x)\|\right)\exp(\Delta L_1).
      \end{aligned}
    \end{equation*}
    This concludes the proof.
  \end{proof}
The rest of the section contains the proofs of the lemmas appearing in the proof of \Cref{thm:complexity}.
Without loss of generality, we may assume that the normalizing constant of the density of $\pi$ is equal to 1. 
Thus, $\pi(x)=e^{-V(x)}$ and  the following lemma is true.

\subsection{Proof of \Cref{lem:upperboundmu0}}
\begin{proof}
  By \Cref{poly-assumption}, we know that  $\big\|\nabla V(x)\big\| \leq \rmQ (\|x\|)$.
Thus, applying Taylor formula
\begin{equation*}\label{eq:aaa}
  \begin{aligned}
    V(x)&=\int_{0}^1\inner{\nabla V(tx)}{x}\rmd t+V(0)\\
    &\leq \int_{0}^1\norm{\nabla V(tx)}\|x\|\rmd t+V(0)\\
    &\leq \int_0^1 \rmQ(\|tx\|) \|x\| \rmd t + V(0).\\
  \end{aligned}
\end{equation*}
Since the coefficients of $\rmQ$ are positive, one may verify that for every $r \geq 0$
\begin{equation*}
  \rmQ(r) \leq \rmQ(1) (r^p + 1).
\end{equation*}
Therefore, 
\begin{equation*}
  \begin{aligned}
    V(x)
     &\leq \int_{0}^1\rmQ(1)\left(\norm{tx}^p+1\right)\|x\|\rmd t+V(0)\\
      &\leq \rmQ(1)\left(\frac{\|x\|^{p+1}}{p+1}+\|x\|\right)+V(0).
  \end{aligned}
\end{equation*}
Now, let us calculate $\KL(\mu_0\mid\pi)$,
  \begin{equation*}
    \begin{aligned}
      \KL(\mu_0\mid\pi)
      &=\int_{\RR^d}\log\left(\frac{\mu_0}{\pi}(x)\right)\mu_0(\rmd x)\\
      &=\int_{\RR^d} \log(\mu_0(x))\mu_0(x)\rmd x+\int_{\RR^d} V(x)\mu_0(\rmd x)\\
      &\leq\frac{d}{2}\log\left(\frac{1}{2\Pi e}\right)+\int_{\RR^d} \left[\rmQ(1)\left(\frac{\|x\|^{p+1}}{p+1}+\|x\|\right)+V(0)\right]\mu_0(\rmd x).
   \end{aligned}
  \end{equation*}

  Combining \Cref{lem:gaussian-moment} with the inequality below
  \begin{equation*}
    \|x\|=\left(\sum_{i=1}^d|x_i|^2\right)^{\frac{1}{2}}
    \leq \left(\sum_{i=1}^d|x_i|\right)^{2\cdot\frac{1}{2}}=\sum_{i=1}^d|x_i|,
    \end{equation*}
   we obtain the following: 
  \begin{equation*}
    \begin{aligned}
      \KL(\mu_0\mid\pi)
      &\leq  \frac{d}{2}\log\left(\frac{1}{2{\Pi} e}\right)+V(0)
      +\frac{\rmQ(1)}{p+1}\cdot \frac{(2d)^{\frac{p+1}{2}}}{\sqrt{ \Pi}}  \Gamma\left(\frac{p+2}{2}\right) 
      + \rmQ(1)\int_{\RR^d} \sum_{i=1}^d|x_i|\mu_0(\rmd x)\\
      &=\frac{d}{2}\log\left(\frac{1}{2{\Pi} e}\right)+V(0)
      +\frac{\rmQ(1)}{p+1}\cdot \frac{(2d)^{\frac{p+1}{2}}}{\sqrt{ \Pi}}  \Gamma\left(\frac{p+2}{2}\right) 
      + \rmQ(1) d \int_{\RR} |r| \frac{1}{\sqrt{2\Pi}}e^{-\frac{r^2}{2}}\rmd r\\
      &=\frac{d}{2}\log\left(\frac{1}{2{\Pi} e}\right)+V(0)
      + \frac{\rmQ(1)(2d)^{\frac{p+1}{2}}  \Gamma\left(\frac{p+2}{2}\right) }{\sqrt{ \Pi}(p+1)} 
      + {\rmQ(1) d}\sqrt{\frac{2}{\Pi}}.
    \end{aligned}
  \end{equation*}
This concludes the proof.
\end{proof}
\subsection{Proof of \Cref{lem:gaussian-moment}}
\begin{proof}
      By Jensen inequality we have 
      \begin{equation*}
        \|{x}\|^{m}=\left(\sum_{i=1}^d|x_i|^2\right)^{\frac{m}{2}}
        =d^{\frac{m}{2}}\left(\frac{1}{d}\sum_{i=1}^d|x_i|^2\right)^{\frac{m}{2}}
        \leq d^{\frac{m-2}{2}}\sum_{i=1}^d|x_i|^{m}.
      \end{equation*}
      Thus, we obtain
      \begin{equation*}
         \begin{aligned}
           \int_{\RR^d} \|{x}\|^{m} \mu_0(\rmd x)
           &\leq  \frac{1}{(2\Pi)^{\nicefrac{d}{2}}}\int_{\RR^d} d^{\frac{m-2}{2}}\sum_{i=1}^d|x_i|^{m} \exp\left(-\frac{\|x\|^2}{2}\right) \rmd x \\
           &\leq  \frac{d^{\frac{m-2}{2}}}{(2\Pi)^{\nicefrac{d}{2}}} \sum_{i=1}^d \int_{\RR^d}|x_i|^{m} 
           \exp\left(-\frac{\|x\|^2}{2}\right) \rmd x \\
           &\leq  \frac{d^{\frac{m-2}{2}}}{\sqrt{2\Pi}} \sum_{i=1}^d \int_{\RR}|x_i|^{m} 
           \exp\left(-\frac{ x_i^2}{2}\right)  \rmd x_i. \\
         \end{aligned}
       \end{equation*} 
      From \cite{winkelbauer2012moments}, we know the $m$-th central absolute moment of the standard one-dimensional Gaussian is equal to ${2^{\nicefrac{m}{2}}\Gamma(\nicefrac{(m+1)}{2})}/{\sqrt{\Pi}}$. 
      Thus, we obtain
      \begin{equation*}
        \begin{aligned}
          \int_{\RR^d} \|{x}\|^{m} \mu_0(x) \rmd x 
          &\leq  \frac{d^{\frac{m-2}{2}}}{\sqrt{2\Pi}} \sum_{i=1}^d \int_{\RR}|x_i|^{m} 
           \exp\left(-\frac{ x_i^2}{2}\right)  \rmd x_i \\
          &\leq  \frac{d^{\frac{m}{2}}}{\sqrt{2\Pi}}  \int_{\RR}|r|^{m} 
           \exp\left(-\frac{ r^2}{2}\right)  \rmd r \\
          &\leq  \frac{(2d)^{\frac{m}{2}}}{\sqrt{ \Pi}}  \Gamma\left(\frac{m+1}{2}\right).\\
        \end{aligned}
      \end{equation*}
  \end{proof}

\section{Miscellaneous}\label{sec:app-assumption}

In this section, we present previously known auxiliary results that were mentioned in the paper. 
Some of the results are proved, the others refer to their papers of origin.

\subsection{Reproducing property}\label{append:repro}
Here we remind you of the reproducing property: let $f\in\cH_0$, then we have $f(x)=\langle f(\cdot),k(x,\cdot)\rangle_{\cH_0}$, if we choose $f(\cdot)=k(y,\cdot)\in\cH_0$, then $k(y,x)=\langle k(y,\cdot),k(x,\cdot)\rangle_{\cH_0}$.

Let us first use the reproducing property to calculate $\norm{k(x,\cdot)}_{\cH_0}$. The definition of the norm in Hilbert spaces yields:
\begin{equation*}
  \norm{k(x,\cdot)}_{\cH_0}^2=\langle k(x,\cdot),k(x,\cdot)\rangle_{\cH_0}=k(x,x).
\end{equation*} 

Now, let us proceed to $\norm{\partial_{x_{i}}k(x,\cdot)}_{\cH_0}$.
Let $e_i$ be the $i$-th standard basis of $\R^d$. Then
\begin{equation*}
  \begin{aligned}
    \langle \partial_{x_{i}}k(x,\cdot)&,\partial_{y_{i}}k(y,\cdot)\rangle_{\cH_0}=\lim_{\epsilon_1\to 0}\lim_{\epsilon_2\to 0}
    \left\langle \frac{k(x+\epsilon_1e_i,\cdot)-k(x,\cdot)}{\epsilon_1}, \frac{k(y+\epsilon_2e_i,\cdot)-k(y,\cdot)}{\epsilon_2}\right\rangle_{\cH_0}\\
    &=\lim_{\epsilon_1\to 0}\lim_{\epsilon_2\to 0}\frac{k(x+\epsilon_1e_i,y+\epsilon_2e_i)-k(x+\epsilon_1e_i,y)-k(x,y+\epsilon_2e_i)+k(x,y)}{\epsilon_1\epsilon_2}\\
    &=\lim_{\epsilon_1\to 0}\lim_{\epsilon_2\to 0}\frac{\left(k(x+\epsilon_1e_i,y+\epsilon_2e_i)-k(x+\epsilon_1e_i,y)\right)-\left(k(x,y+\epsilon_2e_i)-k(x,y)\right)}{\epsilon_1\epsilon_2}\\
    &=\lim_{\epsilon_1\to 0}\frac{\partial_{y_i}k(x+\epsilon_1e_i,y)-\partial_{y_{i}}k(x,y)}{\epsilon_1}\\
    &=\partial_{x_{i}}\partial_{y_{i}}k(x,y).
  \end{aligned}
\end{equation*}
Setting $x=y$,  we obtain $\norm{\partial_{x_{i}}k(x,\cdot)}_{\cH_0}^2=\partial_{x_{i}}\partial_{x_{i}}k(x,x)=:\partial_{x_{i},x_i}^2k(x,x)$, where the first $\partial_{x_{i}}$ is operated on the first variable of $k(\cdot,\cdot)$, the second $\partial_{x_{i}}$ is operated on the second variable of $k(\cdot,\cdot)$.

\subsection{The generalized $\rmT_p$ inequality}\label{sec:app-Tp-assumption}

The following lemma is from \cite{bolley2005weighted}. The proof is omitted.
\begin{lemma}\cite[Corollary 2.3]{bolley2005weighted}\label{bolley}
  Let $\RR^d$ be the Euclidean space with its usual norm. Let $p \geq 1$ and let $\pi$ be a probability measure on $\RR^d$. 
  Assume that there exist $x_{0} \in \RR^d$ and $s>0$ such that $\int_{\RR^d} \exp({s \|x_{0} - x\|^{p}}) \pi(\rmd x)$ is finite. 
  Then
  \begin{equation}\label{eq:bolley-ineq}
    \forall \mu \in \cP(\RR^d), \quad W_{p}(\mu, \pi) 
    \leq \lambda_{BV}\left[\KL(\mu \mid \pi)^{\frac{1}{p}}+\left(\frac{\KL(\mu \mid \pi)}{2}\right)^{\frac{1}{2 p}}\right],
  \end{equation}
  where
  \begin{equation}\label{eq:bolley1}
    \lambda_{BV} :=2 \inf_{x_{0} \in X, s>0}\left(\frac{1}{s}\left(\frac{3}{2}+\log \int_{\RR^d} \exp({s \|x_{0} - x\|^{p}})  \pi(\rmd x)\right)\right)^{\frac{1}{p}}<+\infty.
  \end{equation}
\end{lemma}
 We want to underline the fact that the constant here may depend on the dimension.
Below, we explicit the dimension-dependence of  $\lambda_{BV}$ 
for the general class of distributions  $\pi(x)\propto \exp(-\norm{x}^p)$ with $p\geq 1$. 
   It is straightforward the condition of the lemma is satisfied for $x_0 = 0$.
  Now, let us compute the objective function in \eqref{eq:bolley1}. 
  We start with the integral term:
  \begin{equation*}
    \begin{aligned}
      \int_{x\in\R^d}\exp(s\norm{x}^p)\pi(x)dx
      &=\frac{1}{Z}\int_{x\in\R^d}\exp(-(1-s)\norm{x}^p)dx\\
      &=\frac{1}{Z}(1-s)^{-\frac{d}{p}}\int_{x\in\R^d}\exp(-\norm{x}^p)dx\\
      &=(1-s)^{-\frac{d}{p}}.
    \end{aligned}
  \end{equation*}  
  Inserting the previous formula  into \eqref{eq:bolley1}, we obtain
  \begin{equation*}\label{eq:cpiporder}
    \begin{aligned}
      \lambda_{BV}:&=2 \inf _{x_{0} \in\R^d, s>0}\left(\frac{1}{s}\left(\frac{3}{2}+\log \int \exp({s \norm{x-x_0}^{p}} )d \pi(x)\right)\right)^{\frac{1}{p}}\\
      &=2\inf_{s\in (0,1)}\left(\frac{3}{2s}-\frac{d\log(1-s)}{ps}\right)^{\frac{1}{p}}\\
      &\geq 2\inf_{s\in (0,1)}\left(\frac{3}{2s}+\frac{d}{p}\right)^{\frac{1}{p}}\\
      &=2 d^{\frac{1}{p}}.
    \end{aligned}
  \end{equation*}
  Thus, indeed the constant $\lambda_{BV}$ can be dimension dependent.

\subsection{Grönwall's integral inequality}

The following lemma is the integral form of Gr{\"o}nwall inequality from \cite[Chapter II.]{amann2011ordinary}.

\begin{lemma}[Gr{\"o}nwall Inequality]\label{lem:gronwall}
  Assume $\phi, B:[0, T] \rightarrow \mathbb{R}$ are bounded non-negative measurable function and $C:[0, T] \rightarrow \mathbb{R}$ is a non-negative integrable function with the property that
  \begin{equation}
    \label{eq:grrrr1}
    \phi(t) \leq B(t)+\int_{0}^{t} C(\tau) \phi(\tau) d \tau \quad \text { for all } t \in[0, T]
  \end{equation}
  Then
  \begin{equation}
    \label{eq:grrrr2}
    \phi(t) \leq B(t)+\int_{0}^{t} B(s) C(s) \exp \left(\int_{s}^{t} C(\tau) d \tau\right) d s \quad \text { for all } t \in[0, T].
  \end{equation}
\end{lemma}

\subsection{Stein-Fisher information and weak convergence}\label{sec:app-weak-conv}

  We provide a sufficient condition on which $\lim_{n\to\infty}I_{Stein}(\mu_n\mid\pi)$ implies $\mu_n\to\pi$ weakly. This condition can be found in \cite{gorham2017measuring}. Since Stein Fisher information $I_{Stein}(\cdot\mid\pi)$ depends on the target distribution $\pi$ and the kernel $k(\cdot,\cdot)$, we need the following two properties, respectively:
  
  \begin{enumerate}
    \item $\pi$ is distant dissipative, that is $\kappa_{0} \triangleq \liminf _{r \rightarrow \infty} \kappa(r)>0$  with
    \begin{equation*}
      \kappa(r)=\inf \left\{2 \frac{\langle \nabla V(x)-\nabla V(y), x-y\rangle}{\|x-y\|^{2}}:\|x-y\|=r\right\}.
    \end{equation*}
    If $V$ is strongly convex outside a compact set, then $\pi$ is distant dissipative, for instance $V(x)=\norm{x}^{2+\delta}$ with $\delta\geq 0$.
    \item $k(\cdot,\cdot)$ is an inverse multiquadratic kernel, i.e., $k(x, y)=\left(c^{2}+\|x-y\|^{2}\right)^{\beta}$ for some $c>0$ and $\beta \in(-1,0)$. It is easy to check that \Cref{A-ker-bound} is satisfied.
  \end{enumerate}

\subsection{Proof of \Cref{rem:alt-def-gmu}}\label{sec:app-proof-alt-def-gmu}
Chain rule implies
  \begin{equation*}
    \begin{aligned}
      \int_{\RR^d} \nabla &\log\left(\frac{\mu(x)}{\pi(x)}\right) k(x,\cdot) \mu(\rmd x) \\
      &= \int_{\RR^d} \nabla\log({\mu(x)}) k(x,\cdot) \mu(\rmd x) - \int_{\RR^d} \nabla\log(\pi(x))k(x,\cdot)\mu(\rmd x)\\
      &= \int_{\RR^d}  \nabla {\mu(x)} k(x,\cdot) \rmd x - \int_{\RR^d} \log(\pi(x)) \nabla_x k(x,\cdot)\mu(\rmd x).\\
    \end{aligned}
  \end{equation*}
  Since $\mu(x) \rightarrow 0$, when $\|x\|\rightarrow +\infty$, the integration by parts yields 
   \begin{equation*}
    \begin{aligned}
      \int_{\RR^d}  \nabla {\mu(x)} k(x,\cdot) \rmd x &- \int_{\RR^d} \log(\pi(x)) \nabla_x k(x,\cdot)\mu(\rmd x) \\&=
      - \int_{\RR^d}  \left\{{\mu(x)} \nabla k(x,\cdot) +  \log(\pi(x)) \nabla_x k(x,\cdot)\right\}\mu(\rmd x) \\
      &=g_{\mu}(\cdot).
    \end{aligned}
  \end{equation*}
  This concludes the proof.

\end{document}